\documentclass{article}
\usepackage{amssymb}
\usepackage{amsmath}

 \newcommand{\fiverm}{\tiny }
 \input{prepictex}
 \input{pictex}
 \input{postpictex}

\begin{document}

\title{Smith Normal Form and Acyclic Matrices}
\author{In-Jae Kim and Bryan L. Shader\\ Department of Mathematics\\ University of Wyoming\\ Laramie, WY 82071, USA}

\maketitle

\abstract{An approach, based on  the Smith Normal Form,  is introduced  to study the spectra of symmetric matrices with a given graph.
The  approach serves well to explain how the path cover number (resp. diameter of a tree $T$) is
related to the maximum multiplicity occurring for an eigenvalue of a symmetric matrix whose graph is $T$ (resp. the minimum number $q(T)$ of distinct eigenvalues over the symmetric matrices whose graphs are $T$). The approach is  also  applied to a more general  class of connected graphs $G$, not necessarily trees, in order to establish a lower bound on $q(G)$. }



\font\tenmsb=msbm10
\font\sevenmsb=msbm7
\font\fivemsb=msbm5
\newfam\msbfam
\textfont\msbfam=\tenmsb
\scriptfont\msbfam=\sevenmsb
\scriptscriptfont\msbfam=\fivemsb
\def\Bbb#1{\fam\msbfam\relax#1}
\font\ssfont=cmss10
\newfam\ssfam \def\ss{\fam\ssfam\ssfont} \textfont\ssfam=\ssfont
\def\sh#1{\,{#1}\,}

\newcommand{\Z}{{\Bbb Z}}
\newcommand{\Q}{{\Bbb Q}}
\newcommand{\R}{{\Bbb R}}
\newtheorem{theorem}{Theorem}
\newtheorem{lemma}[theorem]{Lemma}
\newtheorem{corollary}[theorem]{Corollary}
\newtheorem{definition}[theorem]{Definition}
\newtheorem{proposition}[theorem]{Proposition}
\newtheorem{example}[theorem]{Example}
\newtheorem{conjecture}[theorem]{Conjecture}
\newcommand{\mb}{\mathbf}

\section{Introduction}
This paper concerns the relationship between the algebraic and geometric properties
of a symmetric matrix and the combinatorial arrangement of its
nonzero entries (i.e. its graph). We begin by establishing  some basic graph theoretic  notation and terminology that
follows that in  \cite{GR}.

A \emph{graph} $G$ consists of a vertex set $V(G)$ and an \emph{edge} set $E(G)$, where an edge is an unordered
pair of distinct vertices of $G$. We use $uv$ to denote the edge joining vertices $u$ and $v$.
If $uv$ is an edge, then we say that $u$ and $v$ are \emph{adjacent},  and that $v$ is a \emph{neighbor} of $u$.
A vertex is \emph{incident} with an edge if it is one of the two vertices of the edge.
The \emph{degree} of a vertex is the number of  edges incident to the vertex.
A \emph{subgraph} of a graph $G$ is a graph $H$ such that $V(H) \subseteq V(G)$, and $E(H) \subseteq E(G)$.
A subgraph $H$  is an \emph{induced subgraph} if two vertices of $V(H)$ are adjacent in $H$
if and only if they are adjacent in $G$.
If $U \subseteq V(G)$, then $G \setminus U$ denotes the induced subgraph of $G$ whose vertex set is $V(G) \setminus U$.

A \emph{path} $P$ of $G$  is a sequence  $v_{1}, v_2,\ldots, v_{n}$ of
distinct vertices  such that consecutive vertices are adjacent, and is   denoted  by
$v_1$---$v_2$---$\cdots$---$ v_n$. We  say that
$P$ \emph{covers} the vertices $v_{1}, \ldots, v_{n}$, and that $v_{j}$ is \emph{covered by $P$} for  $j=1, \ldots, n$.
The \emph{length} of $P$ is the number of edges in $P$.
If  each vertex of $G$ belongs to at most one of the  paths $P_{1}$, $\ldots$, $P_{k}$, then
$P_{1}, \ldots, P_{k}$ are \emph{disjoint paths}.

If there is a path between each pair of  vertices of   $G$, then $G$ is \emph{connected};
otherwise $G$ is \emph{disconnected}.
A \emph{cycle} is a connected graph where every vertex has exactly two neighbors.
The \emph{length} of a cycle is the number of edges in the cycle.
A \emph{cycle in a graph} $G$ is a subgraph of $G$ that is a cycle.
An \emph{acyclic graph} is a graph with no cycles.
A connected acyclic graph is called a \emph{tree}, and an acyclic graph is called a \emph{forest}.

Let $G$ be a connected graph.
The \emph{distance} between two vertices $u$ and $v$ of $G$ is the  minimum number of edges in
a  path from  $u$ to  $v$.
The \emph{diameter} of a connected graph $G$ is the maximum of the distances over pairs of vertices of $G$, and is denoted by $d(G)$.
If $G$ is a tree, $d(G)$ is the longest length of a  path in $G$.
The \emph{path cover number} of $G$ is the minimum number of disjoint paths needed to cover all of the vertices of $G$, and
is denoted by $p(G)$.

As is customary, we use graphs to model the combinatorial structure of a matrix.
Let $A=[a_{ij}]$ be an $n$ by $n$ symmetric matrix.
The \emph{graph} $G(A)$ of $A$ consists of the vertices $1, 2, \ldots, n$, and  the edges $ij$
for which   $i \ne j$ and $a_{ij} \ne 0$.
Note that $G(A)$ does not depend on the diagonal entries of $A$.
An $n$ by $n$ symmetric matrix $A$ is called an \emph{acyclic} matrix if $G(A)$ is a tree (see \cite{F}).

For a given graph $G$ on $n$ vertices, define $S(G)$ to be the set of all $n$ by $n$ real, symmetric matrices with graph $G$, i.e.
\begin{center}
$S(G)=\{A_{n \times n}$ $|$ $A$ is real, symmetric and $G(A)=G$ $\}$.
\end{center}

For the remainder of this section, matrices are real.
Let $\sigma$ be a multi-list of $n$ real numbers.
If there exists an $n$ by $n$ symmetric matrix $A$ whose spectrum is $\sigma$, then we say that
$\sigma$ is \emph{realized} by $A$, or $A$ \emph{realizes} $\sigma$.
The \emph{spectrum of $S(G)$} for a graph $G$ is the set of all spectra realized by some matrix in $S(G)$.
For a given graph $G$, one can ask to characterize the spectrum of $S(G)$.
This characterization problem is known as the Inverse Eigenvalue Problem
for graphs, or IEP-G for short.
If $G$ is a tree, then we use {IEP-T} instead of IEP-G.

The IEP-G seems quite difficult.
A first step toward  resolving  the IEP-G for a given graph $G$ is to  analyze the
possible multiplicities of the eigenvalues in the spectra of matrices in $S(G)$.
If the distinct eigenvalues of $A$ are $\lambda_{1} < \lambda_{2} < \cdots < \lambda_{q}$ and
their corresponding multiplicities $m_{1}, \ldots, m_{q}$, then $\langle m_{1}, m_{2}, \ldots, m_{q} \rangle$ is
the \emph{ordered multiplicity list} of the eigenvalues of $A$.

The interplay between the spectral properties of acyclic matrices and the combinatorial properties of trees has
been a fruitful area of research for the past 40 years, and some significant, intriguing and recent progress has been made on the IEP-T
(see \cite{F, JD1, JD2, JDS2, JDSSW, LD, LDJ, N}).
The significant graphical parameters of trees $T$ considered for the IEP-T are $d(T)$ and $p(T)$.
For instance, A. Leal Duarte and C.R. Johnson showed in \cite{JD1, LDJ} that the minimum number $q(T)$ of distinct eigenvalues over
the matrices in $S(T)$ satisfies $q(T) \geq d(T)+1$, and the maximum multiplicity $M(T)$ occurring for an eigenvalue of a matrix in $S(T)$
is  $p(T)$.

In this paper, we introduce an approach based on the Smith Normal Form to study the spectra of matrices in $S(G)$. In Section 2 we relate the multiplicities of the eigenvalues of an $n$ by $n$ symmetric matrix $A$
to the Smith Normal Form of $xI-A$ where $I$ is the identity matrix.
In Section 3 we provide a description of the determinant of an $n$ by $n$ matrix in terms of digraphs, and use the description
to show that an eigenvalue of multiplicity $k+1$ or more of an acyclic matrix $A \in S(T)$ for a tree $T$ must be an eigenvalue of each principal submatrix of $A$
whose graph is obtained from $T$ by deleting $k$ disjoint paths.
In addition, as an application, we give an example showing that the IEP-T is not equivalent to determining the ordered multiplicity lists
of the eigenvalues of matrices in $S(T)$.
In Section 4, we show that the tight upper bound $p(T)$ on $M(T)$ is a direct consequence of the Smith Normal Form approach, and
we describe  a systematic way to compute $p(T)$ for a tree $T$.
In Section 5,  the bound $q(T) \geq d(T)+1$ is easily derived, and
 we show that  $q(W) \geq \displaystyle\frac{9d(W)}{8}+\frac{1}{2}$ for an
 infinite family of trees.
In Section 6, we give a lower bound on $q(G)$ for a class of connected graphs $G$.

\section{SNF and Multiplicities of Eigenvalues}
Throughout the remainder of this paper, we let  $I$  denote an
identity matrix of an appropriate order.
In this section we give some useful results on the Smith Normal Form of matrices over the
real polynomial ring, $\mathbb{R}[x]$.
We refer the reader to \cite{BM, HH} for the basic facts.
In particular,  we  relate the multiplicities of the eigenvalues of a real, symmetric matrix
$A$ to the Smith Normal Form of $xI-A$.

For $p(x), q(x) \in \mathbb{R}[x]$ we write {\it $p(x)|q(x)$} if $p(x)$ divides $q(x)$,
and $p(x) \nmid q(x)$ if $p(x)$ does not divide $q(x)$.
We write $(x-a)^k \mid\mid q(x)$ provided $(x-a)^k \mid q(x)$ but $(x-a)^{k+1}\nmid q(x)$.
 We let  $\mathbb{F}$ denote the field of rational functions over $\mathbb{R}$ (that is,
$\mathbb{F}$ is the field of quotients of $\mathbb{R}[x]$), and
let $(\mathbb{R}[x])^{m \times n}$ denote the set of all $m$ by $n$ matrices over $\mathbb{R}[x]$.
Each  $M\in (\mathbb{R}[x])^{m \times n}$ can be viewed as a matrix over $\mathbb{F}$, and
the \emph{rank} of $M$ is defined to be the rank of $M$ over $\mathbb{F}$.
Let $\mbox{GL}_{n}$ be the set of all invertible matrices of order $n$ over $\mathbb{R}[x]$, i.e.
$$
\mbox{\rm{GL}}_{n}= \{ P \in (\mathbb{R}[x])^{n \times n} ~|~ P~ \mbox{has an inverse in} ~(\mathbb{R}[x])^{n \times n} \}.
$$
It is a well-known fact that
\begin{equation}
\label{GL char}
\mbox{GL}_n= \{P\in {\mathbb{R}[x]}^{n\times n}~|~ \det P \in \mathbb{R} \setminus \{ 0 \}\}.
\end{equation}

\noindent
The matrices $M$ and  $N$  in $(\mathbb{R}[x])^{m\times n}$ are \emph{equivalent over $\mathbb{R}[x]$} if there exist $P \in \mbox{GL}_{m}$ and
$Q \in \mbox{GL}_{n}$ such that $N=PMQ$.
Hence equivalent matrices have the same rank.

 A {\it $k$ by $k$ minor of $M$}
is the determinant of a $k$ by $k$ submatrix of $M$.
The monic greatest common divisor of all $k$ by $k$ minors of $M$
is the \emph{$k$th determinantal divisor} of $M$ and is denoted by $\Delta_{k}(M)$.
Another basic fact is:
\begin{proposition}
\label{equivalence}
Let $M, N \in (\mathbb{R}[x])^{n \times n}$.
If $M$ and $N$ are equivalent over $\mathbb{R}[x]$, then
$\Delta_{k}(M) = \Delta_{k}(N)$ for each $k \in \{1, \ldots, n\}$.
\end{proposition}

The following fundamental theorem  asserts that each square matrix over $\mathbb{R}[x]$ is equivalent to a diagonal matrix over $\mathbb{R}[x]$
of a special form (see \cite{BM, HH}).

\begin{theorem}(Smith Normal Form)
\label{SNF}
Let $M \in (\mathbb{R}[x])^{n \times n}$ of rank $r$.
Then there exist $P, Q \in \mbox{\rm{GL}}_{n}$ and monic polynomials $e_{i}(x)$ $(i=1,2,\ldots, r)$ such that $PMQ=D \oplus O$, where
$O$ is the zero matrix of order $n-r$, $D=\mbox{\rm{diag}}(e_{1}(x), \ldots, e_{r}(x))$, and $e_{i}(x)|e_{i+1}(x)$ for $i=1, \ldots, r-1$.
Moreover, $\Delta_{k}(M) = \displaystyle\prod_{j=1}^{k}e_{j}(x)$ and $e_{k}(x) = \frac{\Delta_{k}(M)}{\Delta_{k-1}(M)}$
for each $k \in \{1, \ldots, r\}$.
\end{theorem}

In the above theorem, $D \oplus O$ is called the \emph{Smith Normal Form (SNF)} of $M$, and
$e_{i}(x)$ is called the \emph{$i$th invariant factor} of $M$.

Now assume that $A$ is an $n$ by $n$ real  matrix, and let $S$ be the SNF of $xI-A$.
The characteristic polynomial of $A$, denoted by $p_A(x)$, is $\det (xI-A)=\Delta_n(xI-A)$.
 Since $\det(xI-A)$ is nonzero, the rank of $xI-A$ is $n$.
Thus, $S$ is the full rank matrix of the form $\mbox{diag}(e_{1}(x), \ldots, e_{n}(x))$.
Since $S$ and $xI-A$ are unimodularly equivalent,  by Proposition \ref{equivalence},
\begin{equation}
\label{charac poly}
p_{A}(x)=\displaystyle\prod_{j=1}^{n}e_{j}(x)=\Delta_{n}(S).
\end{equation}

\vskip .1in

Further assume that $A$ is symmetric. Then   the spectrum of a symmetric matrix $A$,
and the  invariant factors of $xI-A$  are closely related.
Let $P, Q \in \mbox{\rm{GL}}_{n}$ such that $P(xI-A)Q=S$.
Since $A$ is symmetric, there exists an orthogonal matrix $U$ of order $n$
so that $U^{T}AU$ is $D=\mbox{diag}(\lambda_{1}, \ldots, \lambda_{n})$.
The diagonal matrix $D$ is called a \emph{diagonalization} of $A$. Thus, $S=PU(xI-D)U^{T}Q$.
Moreover,  since (\ref{GL char}) implies that $PU, U^{T}Q \in \mbox{\rm{GL}}_{n}$, the SNF of $xI-D$ is also $S$.
This along with Proposition~\ref{equivalence} imply that $\Delta_{i}(xI-A) = \Delta_{i}(xI-D) = \Delta_{i}(S)$ for all $i$.

Henceforth $\Delta_{i}(x)$ denotes $\Delta_{i}(xI-A)$. If $\lambda$ is an eigenvalue of $A$, then $m_{A}(\lambda)$ denotes the algebraic multiplicity of $\lambda$.
Let $p(x) \in \mathbb{R}[x]$, and $a \in \mathbb{R}$.
If $(x-a) \mid\mid p(x)$, then $x-a$ is a \emph{linear factor} of $p(x)$.
By considering $xI-D$, we have the following result on the factors of $\Delta_{i}(x)$ and $e_{i}(x)$.

\begin{theorem}
\label{factor main}
Let $A$ be an $n$ by $n$ symmetric matrix whose distinct eigenvalues
are $\lambda_1, \lambda_2, \ldots, \lambda_q$, and  let $S=\mbox{\rm{diag}}(e_{1}(x), \ldots, e_{n}(x))$ be the SNF of $xI-A$.
Then  the following hold:
\begin{itemize}
\item[\rm(a)] if $k\leq n-m_{A}(\lambda_j)$, then  $(x-\lambda_j) \nmid \Delta_{k}(x)$
and $(x-\lambda_j) \nmid e_{k}(x)$,
\item[\rm(b)]
if $k > n-m_{A}(\lambda_j)$, then $(x-\lambda_j)^{k-n+m_{A}(\lambda_j)} || \Delta_{k}(x) $
and   $x-\lambda_j$ is a linear factor of $e_{k}(x)$,
\item[\rm(c)] $e_{n-k}(x)= \prod_{j: m_{A}(\lambda_j)> k} (x-\lambda_j)$
\end{itemize}
\end{theorem}

\noindent
{\bf Proof.}
Fix $j$, and let  $\lambda=\lambda_j$, and $m= m_{A}(\lambda_j)$.
Let $D$ be a diagonalization of $A$.
Without loss of generality, we may assume

$$
xI-D =
\left[ \begin{array}{ccc}
       x-\mu_{1} &         &    \\
                 & \ddots  &    \\
                 &         & x-\mu_{n-m} \end{array} \right]
 \oplus~~~
 \left[ \begin{array}{ccc}
       x-\lambda &         &    \\
                 & \ddots  &    \\
                 &         & x-\lambda \end{array} \right]_{m \times m},
$$
\vskip .1in
\noindent
where  $\mu_{i} \ne \lambda$ for each $i=1, \ldots, n-m$.

\noindent
{\bf (a)}
Suppose $k \leq n-m$. The  determinant of $\mbox{diag}(x-\mu_{1}, \ldots, x-\mu_{k})$
is not divisible by $x-\lambda$.
Thus, $(x-\lambda) \nmid \Delta_{k}(x)$.
By Theorem~\ref{SNF}, $e_{k}(x) | \Delta_{k}(x)$.  Hence,
$(x-\lambda) \nmid e_{k}(x)$.

\noindent
{\bf (b)}
Suppose $k > n-m$, and let $M$ be a $k$ by $k$ submatrix of $xI-D$.
If $M$ is not a principal submatrix of $xI-D$, then $M$ has a zero row and hence $\det M = 0$.
Otherwise, at least $k-(n-m)$ diagonal entries of $M$ are $x-\lambda$.
Thus, $(x- \lambda)^{k-(n-m)} | \Delta_{k}(x)$.
Note that $xI-D$ has a $k$ by $ k$ minor equal to $\det \mbox{diag}(x-\mu_{1}, \ldots, x-\mu_{n-m}, x-\lambda, \ldots, x-\lambda)$.
Thus $(x- \lambda)^{k-(n-m)+1} \nmid \Delta_k(x)$.  Hence
\begin{equation}
\label{exact}
 (x- \lambda)^{k-(n-m)} \mid\mid \Delta_{k}(x).
\end{equation}

By Theorem~\ref{SNF}, $e_{k}(x)=\displaystyle\frac{\Delta_{k}(x)}{\Delta_{k-1}(x)}$.
By (\ref{exact}),  $\Delta_{k}(x)$ has exactly $k-n+m$ factors equal to $x-\lambda$ and $\Delta_{k-1}(x)$ has exactly $(k-n+m)-1$ factors equal to $x-\lambda$.
Thus $(x-\lambda) \mid\mid e_{k}(x)$, and (b) holds.

\noindent
{\bf (c)}  By (\ref{charac poly}),  $e_{n-k}(x)$ is a product of  linear factors
from $\{x-\lambda_{1}, \ldots, x-\lambda_{q}\}$, and by (a) and   (b),
the factors are distinct, and  $x-\lambda_j$ is a factor of $e_{n-k}(x)$ if
and only if $m_{A}(\lambda_{j}) > k$.
  Thus (c) holds.
\hfill{\rule{2mm}{2mm}}

\bigskip\noindent
Useful, immediate consequences  of Theorem  \ref{factor main} are the following:

\begin{corollary}
\label{multi and inv fac}
Let $A$ be an $n$ by $n$ symmetric matrix, and $S=\mbox{\rm{diag}}(e_{1}(x), \ldots, e_{n}(x))$ be the SNF of $xI-A$.
Suppose that $\lambda$ is an eigenvalue of $A$.
Then

\begin{itemize}
\item[\rm(a)]
 $m_{A}(\lambda)\geq k$ if and only if $(x-\lambda) | e_{n-k+1}$
\item[\rm (b)] $m_{A}(\lambda) = k$ if and only if $(x-\lambda) | e_{n-k+1}(x)$ but $(x-\lambda) \nmid e_{n-k}(x)$.
\end{itemize}
\end{corollary}

\noindent
Corollary \ref{multi and inv fac} implies
\begin{equation}
\label{multi k+1}
\mbox{$\mbox{\rm{deg}}(e_{n-k}(x))$ is the number of eigenvalues of $A$ with multiplicity $k+1$ or more.}
\end{equation}

Taking $k=0$ in (c) of Theorem \ref{factor main}, we see that $e_n=(x-\lambda_1)(x-\lambda_2)
\cdots (x-\lambda_q)$, which is known to be the minimal polynomial of $A$.
Thus, $\mbox{\rm deg}(e_n)$ equals the number $q(A)$ of distinct eigenvalues of $A$, and by Theorem \ref{SNF}, we have the following:

\begin{corollary}
\label{min poly}
Let $A$ be an $n$ by $n$ symmetric matrix, and $S=\mbox{\rm{diag}}(e_{1}(x), \ldots, e_{n}(x))$ be the SNF of $xI-A$.
Then $e_{n}(x)$ is the minimal polynomial of $A$, and
$$
q(A)=n-\mbox{\rm{deg}}(\Delta_{n-1}(x)).
$$
\end{corollary}

Corollary~\ref{min poly} allows one to obtain a lower bound on $q(A)$
from an upper bound on $\mbox{deg}(\Delta_{n-1}(x))$.

\section{Eigenvalues of Principal Submatrices of an acyclic matrix}

In this section we associate a digraph on $n$ vertices with an $n$ by $n$ matrix $M$ over $\mathbb{R}[x]$,  and describe $\det M$ in terms of the structure of the digraph associated with $M$.
We use this description to  show that an eigenvalue of a symmetric  acyclic matrix $A$ with multiplicity $k+1$ or more is an eigenvalue
of each principal submatrix of $A$ whose indices correspond to the vertices
not covered by a set of $k$ disjoint paths.
As an application, we provide an example showing that the IEP-T is not equivalent to determining the ordered multiplicity
lists of the eigenvalues of matrices in $S(T)$ (see also \cite{BF}).

First, we give some necessary definitions.
Let $T$ be a tree,  $Q$  an induced subgraph of $T$, and
$M$  a symmetric matrix over $\mathbb{R}[x]$ with $G(M)=T$.
Then $M[Q]$ denotes the principal submatrix of $M$ whose rows and columns correspond to the vertices of $Q$, and  $T \setminus Q$ denotes the induced subgraph of $T$ obtained by deleting all vertices of $Q$. If $Q_{1}, \ldots, Q_{k}$ are induced subgraphs of $T$, $Q_{1} \cup \cdots \cup Q_{k}$ denotes the induced subgraph of $T$
whose vertex set is the union of the vertex sets of $Q_{1}, \ldots, Q_{k}$.

The \emph{end vertices} of the path $v_{1}$---$v_{2}$---$\cdots$--- $v_{\ell}$ are $v_1$ and $v_{\ell}$.
 If  $P_{1} = v_{1}$---$\cdots$ ---$ v_{s}$ and $P_{2} = v_{s}$---$\cdots$---$ v_{t}$ are paths whose
 only common vertex is $v_s$, then $P_{1}P_{2}$ denotes the path $v_{1}$---$ \cdots$--- $v_{s}$--- $\cdots$---$ v_{t}$,
obtained by concatenating $P_1$ and $P_2$.

Let $M=[m_{ij}]$ be an $n$ by $n$ matrix over $\mathbb{R}[x]$.
The \emph{digraph} $D(M)$ of $M$ consists of the vertices $1, 2, \ldots, n$, and arcs $(i,j)$ from vertex $i$ to vertex $j$
if and only if $m_{ij} \ne 0$.
An arc from a vertex to itself is called a \emph{loop}.
A \emph{subdigraph} $H$ of $D(M)$ is a digraph such that the vertex set of $H$ is a subset of $\{1, 2, \ldots, n\}$ and
the arc set of $H$ is a subset of the arc set of $D(M)$.
The \emph{underlying graph} of the digraph $D(M)$ is the graph obtained by treating each arc $(i,j)$ ($i \ne j$) of $D(M)$ as the edge $ij$, and
ignoring the loops.

A \emph{directed walk} in $D(M)$ is a sequence of vertices $(v_{1},v_{2},\ldots, v_{\ell})$,
such that $(v_{i},v_{i+1})$ is an arc for each $i=1, \ldots, \ell-1$, and
$v_{1}$ is the \emph{initial} vertex and $v_{\ell}$ is the \emph{terminal} vertex of the directed walk.
The directed walk $(v_{1}, v_{2},\ldots, v_{\ell})$ \emph{covers} the vertices $v_{1}, v_{2}, \ldots, v_{\ell}$,
and has \emph{length} $\ell-1$.
If $W_{1}$ is the directed walk $(v_{1}, \ldots, v_{s})$, and  $W_{2}$ is the directed walk $(v_{s}, \ldots, v_{t})$, then $(W_{1}, W_{2})$ is the directed walk $(v_{1}, \ldots, v_{s}, \ldots, v_{t})$.
If no vertex of a directed walk is repeated, then the directed walk is a \emph{directed path}.

If the underlying graph of $D(M)$ is a tree $T$, then there is at most one directed path from vertex $i$ to vertex $j$.
If there exists one, $P_{i \to j}$ denotes the unique directed path from vertex $i$ to vertex $j$ in $D(M)$.
Let $P_{i-j}$ denote the unique path connecting $i$ and $j$ in the underlying graph $T$ of $D(M)$.

If the initial vertex is equal to the terminal vertex in a directed walk, then  the directed walk is \emph{closed}.
A \emph{directed cycle} is a closed directed walk with no repeated vertices other than the initial and terminal vertices.
If a directed cycle has length $r$, the directed cycle is a \emph{directed $r$-cycle}.
If each vertex is incident to at most one of the directed cycles $\gamma_{1}, \ldots, \gamma_{k}$ in $D$, then $\gamma_{1}, \ldots, \gamma_{k}$ are \emph{disjoint}.

The \emph{weight of the arc} $(i,j)$ of $D(M)$ is $m_{ij}$.
The \emph{weight of a directed walk} $\beta$ of $D(M)$ is the product of the weights of its arcs, and is denoted by $\mbox{wt}(\beta)$.
The \emph{signed weight of a directed cycle} $\gamma$, denoted by $\mbox{swt}(\gamma)$, is $(-1)^{\ell-1}\mbox{wt}(\gamma)$ where $\ell$ is the length of $\gamma$.
Let $\alpha = \{\gamma_{1}, \ldots, \gamma_{t}\}$ where $\gamma_{1}, \ldots, \gamma_{t}$ are disjoint directed cycles in $D(M)$
covering all of the $n$ vertices.
We define the \emph{signed weight of $\alpha$} to be the product of the signed weights of $\gamma_{1}, \ldots, \gamma_{t}$.
Assume that $\Gamma$ is the set of all such $\alpha$'s described above.
Then, by the definition of the determinant of a square matrix, it is known  \cite{BR} (see p.291) that
\begin{equation}
\label{digraph det}
\det\;M= \displaystyle\sum_{\alpha \in \Gamma}{\mbox{swt}(\alpha)}.
\end{equation}

Let $M$ be an $n$ by $n$ matrix.
If $\alpha$ and $\beta$ are subsets of $ \{1,\ldots,n\}$, then we denote the submatrix of $M$ obtained
by removing (resp. retaining) rows indexed by $\alpha$ and columns indexed by $\beta$ by $M(\alpha, \beta)$ (resp. $M[\alpha, \beta]$).
When $\alpha = \beta$, we use $M(\alpha)$ and $M[\alpha]$, respectively.
We use $e_{i}$ to denote the $i$th column of the identity matrix $I$.

In the following theorem, we provide a relation between
the determinant of a submatrix and the determinant of a principal submatrix of an acyclic matrix over $\mathbb{R}[x]$.

\begin{theorem}
\label{submtrx det}
Let $M$ be an $n$ by $n$ symmetric matrix over $\mathbb{R}[x]$, where the underlying graph of $D(M)$ is a tree $T$, and
let $P_{i_{1}-j_{1}}, \ldots, P_{i_{k}-j_{k}}$ be disjoint paths in $T$.
Then
$$
\det M(\{j_{1}, \ldots, j_{k}\}, \{i_{1}, \ldots, i_{k}\})
=\pm \displaystyle\prod_{s=1}^{k}{\mbox{\rm{wt}}(P_{i_{s} \to j_{s}})} \cdot \det M[T \setminus (P_{i_{1}-j_{1}} \cup \cdots \cup P_{i_{k}-j_{k}})].
$$
\end{theorem}

\noindent
{\bf Proof.}
Let $M'$ be the matrix obtained from $M$ by replacing the $i_{s}$th column of $M$ by $e_{j_{s}}$ for each $s=1, \ldots, k$.
By the construction of $M'$, $M'(\{j_{1}, \ldots, j_{k}\}, \{i_{1}, \ldots, i_{k}\}) = M(\{j_{1}, \ldots, j_{k}\}, \{i_{1}, \ldots, i_{k}\})$.
By  Laplace expansion of the determinant along the columns $i_{1}, i_{2}, \ldots, i_{k}$ of $M'$, we have
\begin{equation}
\label{det M'}
\begin{tabular}{rcl}
$\det M' $ &=&  $\pm \det M'(\{j_{1}, \ldots, j_{k}\}, \{i_{1}, \ldots, i_{k}\})$ \\
           &=&  $\pm \det M(\{j_{1}, \ldots, j_{k}\}, \{i_{1}, \ldots, i_{k}\})$.
\end{tabular}
\end{equation}

In terms of digraphs, $D(M')$ is obtained from $D(M)$ by deleting all incoming arcs to vertices $i_{1}, \ldots, i_{k}$ and
inserting the arcs $(j_{1}, i_{1}), \ldots, (j_{k}, i_{k})$, each with weight 1.
Therefore, $D(M')$ has exactly one arc, namely $(j_{s}, i_{s})$, ending at vertex $i_{s}$ for each $s=1, \ldots, k$.

Set $U=\{(j_{1}, i_{1}), \ldots, (j_{k}, i_{k})\}$.
We claim that each directed cycle in $D(M')$ has at most one arc in $U$.
Suppose to the contrary that there exists a directed cycle $\gamma$ with more than one arc in $U$.
Without loss of generality, we may assume that $(j_{1}, i_{1}), \ldots, (j_{t}, i_{t})$ are arcs of $\gamma$ in $U$, and
$\gamma = ((j_{1}, i_{1}),P_{i_{1} \to j_{2}}, \ldots, (j_{t}, i_{t}), P_{i_{t} \to j_{1}})$ (see Figure 1).

\beginpicture
\small
\setcoordinatesystem units <.90000cm, .90000cm>
\setplotarea x from -2 to 11, y from 1 to 9
\setlinear
\linethickness=0.50pt
\setplotsymbol ({\fiverm .})

\put {$\displaystyle\bullet$} at 3 4
\put {$\displaystyle\bullet$} at 3 7
\put {$\displaystyle\bullet$} at 5 4
\put {$\displaystyle\bullet$} at 5 7
\put {$\displaystyle\bullet$} at 8 4
\put {$\displaystyle\bullet$} at 8 7

\put {$\cdots$} at 6.5 5.5


\plot 3 7 5 4 /
\plot 5 7 6 6 /
\plot 7 5 8 4 /


\plot 2.5 5.5 2.6 5.4 /
\plot 2.5 5.5 2.4 5.4 /

\plot 4.5 5.5 4.6 5.4 /
\plot 4.5 5.5 4.4 5.4 /

\plot 7.5 5.5 7.6 5.4 /
\plot 7.5 5.5 7.4 5.4 /

\plot 4 5.5 4 5.6 /
\plot 4 5.5 3.9 5.5 /

\plot 5.5 6.5 5.5 6.6 /
\plot 5.5 6.5 5.4 6.5 /

\plot 7.5 4.5 7.5 4.6 /
\plot 7.5 4.5 7.4 4.5 /

\plot 9 3 9 3.1 /
\plot 9 3 9.1 3 /

\setquadratic
\plot 3 4  2.5 5.5  3 7 /
\plot 5 4 4.5 5.5 5 7 /
\plot 8 4 7.5 5.5 8 7 /
\plot 3 4 9 3 8 7 /

\put{$i_{1}$} at 3 7.5 
\put{$i_{2}$} at 5 7.5
\put{$i_{t}$} at 8 7.5
\put{$j_{1}$} at 3 3.5
\put{$j_{2}$} at 5 3.5
\put{$j_{t}$} at 8 3.5

\put{$P_{i_{1} \to j_{2}}$} at 3.9 6.5 
\put{$P_{i_{2} \to j_{3}}$} at 6.2 6.5 
\put{$P_{i_{t-1} \to j_{t}}$} at 7 4.2 
\put{$P_{i_{t} \to j_{1}}$} at 9.5 5.5

\put{Figure $1.$} at 6 1.7
\endpicture

\noindent
Since $P_{i_{1} \to j_{2}}, \ldots, P_{i_{t} \to j_{1}}$ are directed paths in the directed cycle $\gamma$,
the underlying graphs of $P_{i_{1} \to j_{2}}, \ldots, P_{i_{t} \to j_{1}}$ are disjoint paths in $T$.

Next, we consider $P_{j_{s} \to i_{s}}$, the unique directed path from $j_{s}$ to $i_{s}$ in $D(M)$ for each $s=1, \ldots, t$.
Then $\tau= (P_{j_{1} \to i_{1}}, P_{i_{1} \to j_{2}},\ldots,  P_{j_{t} \to i_{t}}, P_{i_{t} \to j_{1}})$ is a directed closed walk in $D(M)$.
Since the underlying graph of $D(M)$ is $T$, the multi-set of the arcs in $\tau$ is the disjoint union
of the sets of the arcs of directed 2-cycles.
Note that $P_{j_{1} - i_{1}}, \ldots, P_{j_{t} - i_{t}}$ are disjoint and
that  $P_{i_{1} - j_{2}}, \ldots, P_{i_{t} - j_{1}}$ are disjoint.
Thus, the disjoint union of the sets of the edges of $P_{j_{1}- i_{1}}, \ldots, P_{j_{t}-i_{t}}$ is equal to the disjoint union of the sets of the edges of $P_{i_{1} - j_{2}}, \ldots, P_{i_{t} - j_{1}}$.
Since the  paths $P_{j_{1} - i_{1}}, \ldots, P_{j_{t} - i_{t}}$ in $T$ are  disjoint,
there is no  path from $i_{1}$ to $j_{2}$, consisting of the arcs from $P_{j_{1} \to i_{1}}, \ldots, P_{j_{t} \to i_{t}}$.
This contradicts that $P_{i_{1} \to j_{2}}$ is a directed path from $i_{1}$ to $j_{2}$.
Therefore, each directed cycle in $D(M')$ has at most one arc in $U$.

This implies that  each set of disjoint directed cycles of $D(M')$ that cover every vertex consists of the directed cycles $\beta_{s}=((j_{s}, i_{s}),P_{i_{s} \to j_{s}})$
for all $s=1, \ldots, k$ along with disjoint directed cycles covering every vertex of $T \setminus (P_{i_{1}-j_{1}} \cup \cdots \cup P_{i_{k}-j_{k}})$.
Therefore, (\ref{digraph det}) implies that
$$
\det M' = \displaystyle\prod_{s=1}^{k}{\mbox{swt}(\beta_{s})}  \cdot \displaystyle\sum_{\alpha \in \Gamma(M'[T \setminus (P_{i_{1}-j_{1}} \cup \cdots \cup P_{i_{k}-j_{k}})])}{\mbox{swt}(\alpha)}.
$$
Since $M'[T \setminus (P_{i_{1}-j_{1}} \cup \cdots \cup P_{i_{k}-j_{k}})] = M[T \setminus (P_{i_{1}-j_{1}} \cup \cdots \cup P_{i_{k}-j_{k}})]$,
(\ref{digraph det}) implies
$$
\det M'= \displaystyle\prod_{s=1}^{k}{\mbox{swt}(\beta_{s})} \cdot \det M[T \setminus (P_{i_{1}-j_{1}} \cup \cdots \cup P_{i_{k}-j_{k}})].
$$
Note that $\mbox{swt}(\beta_{s})= \pm \mbox{wt}(P_{i_{s} \to j_{s}})$ for each $s=1, \ldots, k$.
Hence, by (\ref{det M'}),
$$
\begin{tabular}{rcl}
$\det M(\{j_{1}, \ldots, j_{k}\}, \{i_{1}, \ldots, i_{k}\})$ & = & $\pm \det M'$ \\
 & = & $\pm \displaystyle\prod_{s=1}^{k}{\mbox{\rm{wt}}(P_{i_{s} \to j_{s}})} \cdot \det M[T \setminus (P_{i_{1}-j_{1}} \cup \cdots \cup P_{i_{k}-j_{k}})]$.
\end{tabular}
$$
\hfill{\rule{2mm}{2mm}}

\begin{corollary}
\label{main}
Let $A \in S(T)$, where  $T$ is a tree on $n$ vertices, and let $S=\mbox{\rm{diag}}(e_{1}(x), \ldots, e_{n}(x))$ be the SNF of $xI-A$.
If $P_{i_{1}-j_{1}}, \ldots, P_{i_{k}-j_{k}}$ are disjoint paths in $T$ covering $n-t$ vertices of $T$, then
\begin{center}
$\Delta_{n-k}(x) | \det((xI-A)[T \setminus (P_{i_{1}-j_{1}} \cup \cdots \cup P_{i_{k}-j_{k}})])$ and $\mbox{\rm{deg}}(\Delta_{n-k}(x)) \leq t$.
\end{center}
Furthermore, if $\lambda$ is an eigenvalue of $A$ with $m_{A}(\lambda) \geq k+1$, then
$\lambda$ is an eigenvalue of $A[T \setminus (P_{i_{1}-j_{1}} \cup \cdots \cup P_{i_{k}-j_{k}})]$ with multiplicity $m_{A}(\lambda) -k$ or more.
\end{corollary}

\noindent
{\bf Proof.}
Let $P_{i_{s} \to j_{s}}$ be the directed path from $i_{s}$ to $j_{s}$ in $D(xI-A)$ whose underlying graph is $P_{i_{s}-j_{s}}$ for each $s=1, \ldots, k$.
By Theorem~\ref{submtrx det},
$$
\det[(xI-A)(\{j_{1}, \ldots, j_{k}\}, \{i_{1}, \ldots, i_{k}\})]
=\pm \displaystyle\prod_{s=1}^{k}{\mbox{\rm{wt}}(P_{i_{s} \to j_{s}})} \cdot \det((xI-A)[T \setminus (P_{i_{1}-j_{1}} \cup \cdots \cup P_{i_{k}-j_{k}})]).
$$
Note that $\mbox{wt}(P_{i_{s} \to j_{s}})$ is a nonzero constant for each $s=1, \ldots, k$.
Hence,
$$
\det[(xI-A)(\{j_{l}, \ldots, j_{k} \}, \{i_{1}, \ldots, i_{k} \})] = c \cdot \det((xI-A)[T \setminus (P_{i_{1}-j_{1}} \cup \cdots \cup P_{i_{k}-j_{k}})])
$$
for some nonzero constant $c$.

Since $(xI-A)(\{j_{1}, \ldots, j_{k} \}, \{i_{1}, \ldots, i_{k} \})$
is an $n-k$ by $n-k$ submatrix of $xI-A$,
\begin{equation}
\label{Delta n-k}
\Delta_{n-k}(x) | \det((xI-A)[T \setminus (P_{i_{1}-j_{1}} \cup \cdots \cup P_{i_{k}-j_{k}})]).
\end{equation}
By Theorem~\ref{factor main}, if $m_{A}(\lambda) \geq k+1$, then
the multiplicity of $(x-\lambda)$ as a factor of $\Delta_{n-k}(x)$ is $m_{A}(\lambda) - k$.
Therefore, by (\ref{Delta n-k}), $(x - \lambda)^{m_{A}(\lambda) - k} | \det((xI-A)[T \setminus (P_{i_{1}-j_{1}} \cup \cdots \cup P_{i_{k}-j_{k}})])$.
This implies that $\lambda$ is an eigenvalue of $A[T \setminus (P_{i_{1}-j_{1}} \cup \cdots \cup P_{i_{k}-j_{k}})]$ with multiplicity $m_{A}(\lambda) -k$ or more.

Note that the degree of $\det((xI-A)[T \setminus (P_{i_{1}-j_{1}} \cup \cdots \cup P_{i_{k}-j_{k}})])$ is $t$.
Thus, by (\ref{Delta n-k}),
$$
\mbox{\rm{deg}}(\Delta_{n-k}(x)) \leq t.
$$
\hfill{\rule{2mm}{2mm}}

\bigskip\noindent
If $\mbox{\rm{deg}}(\Delta_{n-k}(x)) \leq t$, then Theorem~\ref{SNF} implies $\mbox{\rm{deg}}(e_{n-k}(x)) \leq t$.
Hence, Corollary~\ref{multi k+1} implies the following.

\begin{corollary}
\label{main cor}
Let $A \in S(T)$, where  $T$ is a tree on $n$ vertices.
If $k$ disjoint paths of $T$ cover $n-t$ vertices of $T$, then
there are at most $t$ eigenvalues of $A$ with multiplicity $k+1$ or more.
\end{corollary}

It was conjectured in \cite{JDS2} that the IEP-T for a tree $T$ is equivalent to determining
the ordered multiplicity lists of the eigenvalues of matrices in $S(T)$, i.e.~each multi-list of real numbers having an ordered multiplicity list
of the eigenvalues of a matrix in $S(T)$ is the spectrum of a matrix in $S(T)$.
Indeed, \cite{JDS2} showed that for some classes of trees, these two problems are equivalent.
A counterexample to the conjecture was given in  \cite{BF}.  We give a counterexample
on fewer vertices, and a simple argument that shows this is a counterexample.

Consider the tree $T$ illustrated in Figure 2.
We will  show that an ordered multiplicity list of the eigenvalues of a matrix in $S(T)$
requires the eigenvalues having the ordered multiplicity list
to satisfy a certain algebraic condition.

\begin{example}
\label{shaun}

\beginpicture
\small
\setcoordinatesystem units <.80000cm, .80000cm>
\setplotarea x from -4 to 10, y from -.5 to 8.5
\setlinear
\linethickness=0.50pt
\setplotsymbol ({\fiverm .})

\put {$\displaystyle\bullet$} at 3.5 7
\put {$\displaystyle\bullet$} at 6.5 7
\put {$\displaystyle\bullet$} at 5 6
\put {$\displaystyle\bullet$} at 5 4.5
\put {$\displaystyle\bullet$} at 1.5 3.5

\put {$\displaystyle\bullet$} at 3.5 3.5
\put {$\displaystyle\bullet$} at 6.5 3.5
\put {$\displaystyle\bullet$} at 8.5 3.5
\put {$\displaystyle\bullet$} at 3.5 1.5
\put {$\displaystyle\bullet$} at 6.5 1.5

\plot 3.5 7 5 6 /
\plot 5 6 6.5 7 /
\plot 5 6 5 4.5 /

\plot 5 4.5 3.5 3.5 /
\plot 1.5 3.5 3.5 3.5 /
\plot 3.5 3.5 3.5 1.5 /

\plot 5 4.5 6.5 3.5 /
\plot 6.5 3.5 8.5 3.5 /
\plot 6.5 3.5 6.5 1.5 /

\put{$4$} at 3 7
\put{$5$} at 7 7
\put{$1$} at 5.5 5.8
\put{$6$} at 5.5 4.5
\put{$7$} at 1 3.5
\put{$2$} at 3.5 3.9
\put{$3$} at 6.5 3.9
\put{$10$} at 9 3.5
\put{$8$} at 3 1.5
\put{$9$} at 7 1.5
\put{\rm{Figure} $2.$} at 5 .3
\endpicture 

It can be verified that the eigenvalues of the following matrix $A$ in $S(T)$ are $-\sqrt{5}$, $-\sqrt{2}$, $0$, $\sqrt{2}$, $\sqrt{5}$,
and the ordered multiplicity list of the eigenvalues of $A$ is $\langle 1, 2, 4, 2,1 \rangle$:

$$
\left[ \begin{array}{ccccccccccc}
0 & 0 & 0 & \vline & 1 & 1 & 1 & 0 & 0 & 0 & 0 \\
0 & 0 & 0 & \vline & 0 & 0 & 1 & 1 & 1 & 0 & 0 \\
0 & 0 & 0 & \vline & 0 & 0 & 1 & 0 & 0 & 1 & 1 \\
\noalign{\hrule height 0.3pt}
1 & 0 & 0 & \vline & 0 & 0 & 0 & 0 & 0 & 0 & 0 \\
1 & 0 & 0 & \vline & 0 & 0 & 0 & 0 & 0 & 0 & 0 \\
1 & 1 & 1 & \vline & 0 & 0 & 0 & 0 & 0 & 0 & 0 \\
0 & 1 & 0 & \vline & 0 & 0 & 0 & 0 & 0 & 0 & 0 \\
0 & 1 & 0 & \vline & 0 & 0 & 0 & 0 & 0 & 0 & 0 \\
0 & 0 & 1 & \vline & 0 & 0 & 0 & 0 & 0 & 0 & 0 \\
0 & 0 & 1 & \vline & 0 & 0 & 0 & 0 & 0 & 0 & 0 \\
\end{array} \right]$$
\vskip .2in

Suppose that
$$
\sigma =(\lambda_{1},\lambda_{2},\lambda_{2},\lambda_{3},\lambda_{3},\lambda_{3},\lambda_{3},\lambda_{4},\lambda_{4},\lambda_{5})
$$
is realized by a matrix $B$ in $S(T)$ as its spectrum where $\lambda_{1} < \lambda_{2} < \lambda_{3} < \lambda_{4} < \lambda_{5}$.

The disjoint paths $4$---$1$---$5$, $7$---$2$---$8$; and $10$---$3$---$9$ cover all the vertices except vertex $6$ (see Figure $3$).

\beginpicture
\small
\setcoordinatesystem units <.80000cm, .80000cm>
\setplotarea x from -4 to 10, y from -.5 to 8
\setlinear
\linethickness=0.50pt
\setplotsymbol ({\fiverm .})

\put {$\displaystyle\bullet$} at 3.5 7
\put {$\displaystyle\bullet$} at 6.5 7
\put {$\displaystyle\bullet$} at 5 6
\put {$\displaystyle\bullet$} at 5 4.5
\put {$\displaystyle\bullet$} at 1.5 3.5

\put {$\displaystyle\bullet$} at 3.5 3.5
\put {$\displaystyle\bullet$} at 6.5 3.5
\put {$\displaystyle\bullet$} at 8.5 3.5
\put {$\displaystyle\bullet$} at 3.5 1.5
\put {$\displaystyle\bullet$} at 6.5 1.5

\plot 3.5 7 5 6 /
\plot 5 6 6.5 7 /

\plot 1.5 3.5 3.5 3.5 /
\plot 3.5 3.5 3.5 1.5 /

\plot 6.5 3.5 8.5 3.5 /
\plot 6.5 3.5 6.5 1.5 /

\put{$4$} at 3 7
\put{$5$} at 7 7
\put{$1$} at 5.5 5.8
\put{$6$} at 5.5 4.5
\put{$7$} at 1 3.5
\put{$2$} at 3.5 3.9
\put{$3$} at 6.5 3.9
\put{$10$} at 9 3.5
\put{$8$} at 3 1.5
\put{$9$} at 7 1.5
\put{\rm{Figure} $3.$} at 5 .3
\endpicture

\noindent
Since $m_{B}(\lambda_{3}) =4$, Corollary~$\ref{main}$ implies that $\lambda_{3}$ is the eigenvalue of $B[6]$, i.e. $B[6] = \lambda_{3}$.

Next, since the three disjoint paths $4$---$1$---$6$---$2$---$7$; $10$---$3$---$9$; and $8$ cover all the vertices except vertex $5$ (see Figure $4$),
Corollary~$\ref{main}$ implies that $\lambda_{3}$ is the eigenvalue of $B[5]$, i.e. $B[5] = \lambda_{3}$.

\beginpicture
\small
\setcoordinatesystem units <.80000cm, .80000cm>
\setplotarea x from -4 to 10, y from -.5 to 8.5
\setlinear
\linethickness=0.50pt
\setplotsymbol ({\fiverm .})

\put {$\displaystyle\bullet$} at 3.5 7
\put {$\displaystyle\bullet$} at 6.5 7
\put {$\displaystyle\bullet$} at 5 6
\put {$\displaystyle\bullet$} at 5 4.5
\put {$\displaystyle\bullet$} at 1.5 3.5

\put {$\displaystyle\bullet$} at 3.5 3.5
\put {$\displaystyle\bullet$} at 6.5 3.5
\put {$\displaystyle\bullet$} at 8.5 3.5
\put {$\displaystyle\bullet$} at 3.5 1.5
\put {$\displaystyle\bullet$} at 6.5 1.5

\plot 3.5 7 5 6 /
\plot 5 6 5 4.5 /

\plot 5 4.5 3.5 3.5 /
\plot 1.5 3.5 3.5 3.5 /

\plot 6.5 3.5 8.5 3.5 /
\plot 6.5 3.5 6.5 1.5 /

\put{$4$} at 3 7
\put{$5$} at 7 7
\put{$1$} at 5.5 5.8
\put{$6$} at 5.5 4.5
\put{$7$} at 1 3.5
\put{$2$} at 3.5 3.9
\put{$3$} at 6.5 3.9
\put{$10$} at 9 3.5
\put{$8$} at 3 1.5
\put{$9$} at 7 1.5

\put{\rm{Figure} $4.$} at 5 .3
\endpicture

\noindent
Similarly, $\lambda_{3}$ is the eigenvalue of $B[i]$, i.e. $B[i] = \lambda_{3}$ for each $i=4, 7, 8, 9, 10$.

The three disjoint  paths $7$---$2$---$6$---$3$---$10$; $8$ and $9$ cover seven  vertices (see Figure $5$).

\beginpicture
\small
\setcoordinatesystem units <.80000cm, .80000cm>
\setplotarea x from -4 to 10, y from -.5 to 8
\setlinear
\linethickness=0.50pt
\setplotsymbol ({\fiverm .})

\put {$\displaystyle\bullet$} at 3.5 7
\put {$\displaystyle\bullet$} at 6.5 7
\put {$\displaystyle\bullet$} at 5 6
\put {$\displaystyle\bullet$} at 5 4.5
\put {$\displaystyle\bullet$} at 1.5 3.5

\put {$\displaystyle\bullet$} at 3.5 3.5
\put {$\displaystyle\bullet$} at 6.5 3.5
\put {$\displaystyle\bullet$} at 8.5 3.5
\put {$\displaystyle\bullet$} at 3.5 1.5
\put {$\displaystyle\bullet$} at 6.5 1.5

\plot 3.5 7 5 6 /
\plot 5 6 6.5 7 /

\plot 5 4.5 3.5 3.5 /
\putrule from 1.5 3.5 to 3.5 3.5

\plot 5 4.5 6.5 3.5 /
\putrule from 6.5 3.5 to 8.5 3.5

\put{$4$} at 3 7
\put{$5$} at 7 7
\put{$1$} at 5.5 5.8
\put{$6$} at 5.5 4.5
\put{$7$} at 1 3.5
\put{$2$} at 3.5 3.9
\put{$3$} at 6.5 3.9
\put{$10$} at 9 3.5
\put{$8$} at 3 1.5
\put{$9$} at 7 1.5

\put{\rm{Figure} $5.$} at 5 .3
\endpicture

\noindent
Thus, by Corollary~$\ref{main}$, $\lambda_{3}$ is an eigenvalue of $B[\{4,1,5\}]$.
Similarly, $\lambda_{3}$ is an eigenvalue of $B[\{10,3,9\}]$ and $B[\{7,2,8\}]$.

Finally, we show that both eigenvalues of $B$ with multiplicity $2$ are eigenvalues of $B[\{4,1,5\}]$, $B[\{7,2,8\}]$ and $B[\{10,3,9\}]$.
We consider only the case for $B[\{4,1,5\}]$.
Since the single path $7$---$2$---$6$---$3$---$10$ covers five vertices, and $m_{B}(\lambda_{2}) = m_{B}(\lambda_{4})=2$,
Corollary~$\ref{main}$ implies that $\lambda_{2}$ and $\lambda_{4}$ are eigenvalues of $B[\{4,1,5,8,9\}]$.
However, $B[\{4,1,5,8,9\}]= B[\{4,1,5\}] \oplus B[8] \oplus B[9]$.
Since $B[8]=B[9]=\lambda_{3}$, $\lambda_{2}$ and $\lambda_{4}$ are eigenvalues of $B[\{4,1,5\}]$.
The other cases can be shown by choosing the paths $4$---$1$---$6$---$2$---$7$ and $4$---$1$---$6$---$3$---$10$, respectively.

So far, we have shown that $\sigma(B[6])=\lambda_{3}$, and
$\sigma(B[\{4,1,5\}])=\sigma(B[\{7,2,8\}])=\sigma(B[\{10,3,9\}])=(\lambda_{2}, \lambda_{3}, \lambda_{4})$.
Now, we consider the trace of $B$.
The trace of $B$ is equal to the sum of the traces of $B[6]$, $B[\{4,1,5\}]$, $B[\{7,2,8\}]$ and $B[\{10,3,9\}]$.
Since the trace is equal to the sum of all the eigenvalues, we have
$$
\lambda_{1} + 2\lambda_{2} + 4 \lambda_{3} + 2\lambda_{4} + \lambda_{5} = \lambda_{3} + 3 (\lambda_{2} + \lambda_{3} + \lambda_{4})
$$
and hence,
$$
\lambda_{1} + \lambda_{5} = \lambda_{2} + \lambda_{4}.
$$

Therefore, if the ordered multiplicity list $\langle 1,2,4,2,1 \rangle$ is realized by a matrix $B$ in $S(T)$, then
$\sigma(B) = (\lambda_{1},\lambda_{2},\lambda_{2},\lambda_{3},\lambda_{3},\lambda_{3},\lambda_{3},\lambda_{4},\lambda_{4},\lambda_{5})$
must satisfy $\lambda_{1} + \lambda_{5} = \lambda_{2} + \lambda_{4}$.
For instance, $\sigma=(2,3,3,5,5,5,5,7,7,10)$ with the ordered multiplicity list $\langle 1,2,4,2,1 \rangle$
cannot be realized by any matrix in $S(T)$.
\end{example}

\section{Maximum Multiplicity and $p(T)$}
Let $T$ be a tree on $n$ vertices, and let $M(G)$ to denote the maximum multiplicity occurring for an eigenvalue among the  matrices in $S(T)$. Recall that the path cover number of $T$ is  the minimum number of disjoint paths that cover $T$.

Now let $A \in S(T)$.
If $k$ disjoint paths in $T$ cover all the vertices of $T$,
then Corollary~\ref{main cor} implies that no eigenvalue of $A$ has multiplicity more than $k$.
Since $A$ is an arbitrary matrix in $S(T)$, we have the following corollary.

\begin{corollary}
\label{bound for M}
Let $T$ be a tree.
Then
$$
M(T) \leq p(T).
$$
\end{corollary}

Furthermore, it was shown in \cite{JD1} that the upper bound is tight.

\begin{theorem}
\label{M=p}
Let $T$ be a tree.
Then
$$
M(T) = p(T).
$$
\end{theorem}

For the remainder of this section, we describe a systematic way of computing $p(T)$ for a tree $T$.
This method will be used repeatedly in   the following sections.
We first show that the existence of a specific path
for a given tree (see also \cite[Lemma 3.1]{N}).

\begin{proposition}
\label{special path}
Let $T$ be a tree on $n$ vertices.
Then there exists a path in $T$ such that the end vertices of the path are pendant vertices of $T$, and
at most one vertex of the path has degree $3$ or more in $T$.
\end{proposition}

\noindent
{\bf Proof.}
The proof is by induction on $n$.  The result is clear if $n\leq 2$. Assume $n\geq 3$
and proceed by induction.

If $T$ has diameter $2$, then any path of length $2$ works.
Assume that the diameter of $T$ is at least 3, and
let $P$ be a path in $T$, $u$---$v$---$a_1$---$\cdots$---$a_{k-1}$---$a_k$, whose length is the diameter of $T$.
Since $u$ is a pendant vertex of $T$, $T \setminus \{u\}$ is also a tree.
By the inductive hypothesis,  there exists a path $P'$ of $T \setminus \{u\}$ satisfying the given condition.

If the path $P'$ in $T \setminus \{u\}$ does not contain $v$, then $P'$ is also a path in $T$ satisfying the given condition.

If $P'$ contains $v$ as an end vertex, then the path $P_{u-v}P'$ is a path in $T$ satisfying the given condition.

Otherwise $P'$ contains $v$ and $v$ is not an end vertex of $P'$.  Thus there exists a neighbor $w$ of $v$ in $T \setminus \{u\}$
such that $w \ne a_1$.
Suppose that $w$ is not a pendant vertex of $T$.
Then there exists a neighbor $y$ of $w$ other than $v$.
Since $T$ is a tree $y\notin\{a_2, \ldots, a_k\}$ and hence,
$y$---$w$---$v$---$a_1$---$\cdots$---$a_{k-1}$---$a_k$ is a path in $T$, whose length is longer than $P$, which is a contradiction.
Therefore, $w$ is a pendant vertex of $T$, and the path $u$---$v$---$w$ in $T$ satisfies the given condition.
\hfill{\rule{2mm}{2mm}}
\vskip .1in

If $T$ is a path, then a path in $T$ satisfying the conditions in Proposition~\ref{special path} is $T$ itself.
Otherwise, there exists a path $P$ in $T$ such that the end vertices of $P$ are pendant vertices of $T$, and
exactly one vertex of $P$ has degree $3$ or more in $T$.
Next, we show that for such a path $P$, $p(T \setminus P) = p(T) -1$.

\begin{proposition}
\label{path cover number}
Let $T$ be a tree which is not a path.
Suppose that $P$ is a path in $T$ such that $P$'s end vertices are pendant vertices of $T$ and
$P$ has exactly one vertex $v$ of degree $3$ or more in $T$.
Then
$$
p(T)= p(T \setminus P) +1.
$$
\end{proposition}

\noindent
{\bf Proof.}
Note that each  path cover of $T\setminus P$ can be extended to a path cover of $T$
by including the path $P$. Hence $p(T) \leq p(T \setminus P) + 1$

Now we show $p(T) \geq p(T \setminus P) + 1$.
Let $p$ denote $p(T)$. and $C=\{P_{i}\}_{i=1}^{p}$ be a set of $p$ disjoint paths in $T$ covering all of the vertices of $T$.

If $P \in C$,
then, since $C \setminus \{P\}$ covers all of the vertices of $T \setminus P$, $p(T \setminus P) \leq p-1$.

Otherwise, $P \not\in C$.
Then two disjoint paths in $C$, say $\alpha$ and $\beta$, are needed to cover the vertices of $P$.
Assume that $\beta$ covers the vertex $v$.
Then $\alpha$ covers only the vertices of $P$.
Thus, $(C \setminus \{\alpha, \beta\}) \cup \{\beta \setminus P \}$ covers all the vertices of $T \setminus P$.
This implies that $p(T \setminus P) \leq p-1$; equivalently,
$
p \geq p(T \setminus P) + 1.
$
\hfill{\rule{2mm}{2mm}}

\medskip
\noindent
By repeated use of  Proposition~\ref{path cover number},  we can effectively compute the path cover number of the tree $T$ in Figure 6.

\begin{example}
\label{computing p}

\beginpicture
\small
\setcoordinatesystem units <.800000cm, .800000cm>
\setplotarea x from -4 to 10, y from -.5 to 5
\setlinear
\linethickness=0.50pt
\setplotsymbol ({\fiverm .})

\put {$\displaystyle\bullet$} at 1 4.5
\put {$\displaystyle\bullet$} at 2 3
\put {$\displaystyle\bullet$} at 1 1.5
\put {$\displaystyle\bullet$} at 4 3
\put {$\displaystyle\bullet$} at 6 3

\put {$\displaystyle\bullet$} at 4 1.5
\put {$\displaystyle\bullet$} at 6 1.5
\put {$\displaystyle\bullet$} at 9 4.5
\put {$\displaystyle\bullet$} at 8 3
\put {$\displaystyle\bullet$} at 9 1.5

\plot 1 4.5 2 3 /
\plot 2 3 1 1.5 /

\plot 2 3 4 3 /
\plot 4 3 4 1.5 /
\plot 4 3 6 3 /
\plot 6 3 6 1.5 /
\plot 6 3 8 3 /

\plot 8 3 9 4.5 /
\plot 8 3 9 1.5 /

\put{$3$} at .5 4.5
\put{$1$} at 1.5 3
\put{$4$} at .5 1.5
\put{$6$} at 4 3.5
\put{$7$} at 6 3.5
\put{$5$} at 3.5 1.5
\put{$8$} at 6.5 1.5
\put{$9$} at 9.5 4.5
\put{$2$} at 8.5 3
\put{$10$} at 9.5 1.5

\put{\rm{Figure} $6.$} at 5 .3
\endpicture 

\noindent
Since $P_{1}$, $3$---$1$---$4$, and $P_{2}$, $9$---$2$---$10$, satisfy the condition in Proposition~$\ref{path cover number}$ for $T$ and $T \setminus P_{1}$,
respectively, we have the following three disjoint paths covering all of the vertices of $T$.

\beginpicture
\small
\setcoordinatesystem units <.800000cm, .800000cm>
\setplotarea x from -4 to 10, y from -.5 to 5.5
\setlinear
\linethickness=0.50pt
\setplotsymbol ({\fiverm .})

\put {$\displaystyle\bullet$} at 1 4.5
\put {$\displaystyle\bullet$} at 2 3
\put {$\displaystyle\bullet$} at 1 1.5
\put {$\displaystyle\bullet$} at 4 3
\put {$\displaystyle\bullet$} at 6 3

\put {$\displaystyle\bullet$} at 4 1.5
\put {$\displaystyle\bullet$} at 6 1.5
\put {$\displaystyle\bullet$} at 9 4.5
\put {$\displaystyle\bullet$} at 8 3
\put {$\displaystyle\bullet$} at 9 1.5

\plot 1 4.5 2 3 /
\plot 2 3 1 1.5 /

\plot 4 3 4 1.5 /
\plot 4 3 6 3 /
\plot 6 3 6 1.5 /

\plot 8 3 9 4.5 /
\plot 8 3 9 1.5 /

\put{$3$} at .5 4.5
\put{$1$} at 1.5 3
\put{$4$} at .5 1.5
\put{$6$} at 4 3.5
\put{$7$} at 6 3.5
\put{$5$} at 3.5 1.5
\put{$8$} at 6.5 1.5
\put{$9$} at 9.5 4.5
\put{$2$} at 8.5 3
\put{$10$} at 9.5 1.5

\put{\rm{Figure} $7.$} at 5 .3
\endpicture 

\noindent
Hence, $p(T) =3$.
This implies, by Theorem~$\ref{M=p}$, that if a multi-list $\sigma$ of $10$ real numbers has an element with multiplicity greater than $3$,
$\sigma$ cannot be realized by any matrix in $S(T)$.
\end{example}

\section{Relationship between $d(T)$ and $q(T)$}
Let $T$ be a tree, and $A \in S(T)$.
The number of distinct eigenvalues of $A$ is denoted by $q(A)$ and, $q(T)$ denotes the minimum of $q(A)$ over all $A \in S(T)$, and
$d(T)$ is the diameter of $T$.
In this section, we study the relation between $d(T)$ and $q(T)$.

Let $\Delta_{i}(x)$ be the $i$th determinantal divisor of $xI-A$.
By Corollary~\ref{min poly},
\begin{equation}
\label{deg}
q(A) = n - \mbox{deg}(\Delta_{n-1}(x))
\end{equation}
 and
by Corollary~\ref{main}, if a path in $T$ has $\ell$ vertices, then
\begin{equation}
\label{deg Delta n-1}
\mbox{deg}(\Delta_{n-1}(x)) \leq n-\ell.
\end{equation}
Thus, by choosing a path of the longest length, (\ref{deg}) and (\ref{deg Delta n-1}) imply that
\begin{equation}
\label{q(A)}
q(A) \geq d(T)+1.
\end{equation}
Since (\ref{q(A)}) holds for every matrix in $S(T)$, $q(T) \geq d(T)+1$.
Thus,  the following  known theorem (see  \cite{LDJ}) follows easily
from the Smith Normal Form approach.

\begin{theorem}
\label{q(T)}
Let $T$ be a tree.
Then
$$
q(T) \geq d(T)+1.
$$
\end{theorem}

Next, we provide a class of trees $W$ for which $q(W)$ is much larger than $d(W)+1$.
The \emph{$(3,\ell)$-whirl} ($\ell \geq 1$), $W$, is the tree on $n=6\ell+4$ vertices with 6 pendant paths $\alpha_{i}, \beta_{i}, \gamma_{i}$,
each with $\ell$ vertices, for $i=1,2$ as illustrated in Figure 8.
      The vertex $v$ in Figure 8 is the \emph{axis vertex} of $W$, and each of the pendant paths is a \emph{leg} of $W$.

\beginpicture
\small
\setcoordinatesystem units <.700000cm, .700000cm>
\setplotarea x from -4.5 to 12, y from -2 to 10.5
\setlinear
\linethickness=0.50pt
\setplotsymbol ({\fiverm .})

\put {$\displaystyle\bullet$} at 2 8
\put {$\displaystyle\bullet$} at 2.5 7.5
\put {$\displaystyle\bullet$} at 3.5 6.5
\put {$\displaystyle\bullet$} at 4 6

\put {$\displaystyle\bullet$} at 3 9
\put {$\displaystyle\bullet$} at 3.5 8.5
\put {$\displaystyle\bullet$} at 4.5 7.5
\put {$\displaystyle\bullet$} at 5 7

\put {$\displaystyle\bullet$} at 5 6
\put {$\displaystyle\bullet$} at 6 4

\put {$\displaystyle\bullet$} at 9 9
\put {$\displaystyle\bullet$} at 8.5 8.5
\put {$\displaystyle\bullet$} at 7.5 7.5
\put {$\displaystyle\bullet$} at 7 7

\put {$\displaystyle\bullet$} at 10 8
\put {$\displaystyle\bullet$} at 9.5 7.5
\put {$\displaystyle\bullet$} at 8.5 6.5
\put {$\displaystyle\bullet$} at 8 6

\put {$\displaystyle\bullet$} at 7 6

\put {$\displaystyle\bullet$} at 6.7 3.3
\put {$\displaystyle\bullet$} at 6.7 2.5
\put {$\displaystyle\bullet$} at 6.7 1.5
\put {$\displaystyle\bullet$} at 6.7 .7

\put {$\displaystyle\bullet$} at 5.3 3.3
\put {$\displaystyle\bullet$} at 5.3 2.5
\put {$\displaystyle\bullet$} at 5.3 1.5
\put {$\displaystyle\bullet$} at 5.3 .7

\put {$\displaystyle\bullet$} at 6 5

\plot 2 8 2.5 7.5 /
\plot 3.5 6.5 4 6 /
\plot 4 6 5 6 /

\plot 3 9 3.5 8.5 /
\plot 4.5 7.5 5 7 /
\plot 5 7 5 6 /

\plot 9 9 8.5 8.5 /
\plot 7.5 7.5 7 7 /
\plot 7 7 7 6 /

\plot 10 8 9.5 7.5 /
\plot 8.5 6.5 8 6 /
\plot 8 6 7 6 /

\plot 5 6  6 5 /
\plot 6 5 7 6 /
\plot 6 5 6 4 /

\plot 6 4 5.3 3.3 /
\plot 5.3 3.3 5.3 2.5 /
\plot 5.3 1.5 5.3 .7 /

\plot 6 4 6.7 3.3 /
\plot 6.7 3.3 6.7 2.5 /
\plot 6.7 1.5 6.7 .7 /

\put{$\cdot$} at 3 7
\put{$\cdot$} at 2.9 7.1
\put{$\cdot$} at 3.1 6.9

\put{$\cdot$} at 4 8
\put{$\cdot$} at 3.9 8.1
\put{$\cdot$} at 4.1 7.9

\put{$\cdot$} at 8 8
\put{$\cdot$} at 8.1 8.1
\put{$\cdot$} at 7.9 7.9

\put{$\cdot$} at 9 7
\put{$\cdot$} at 9.1 7.1
\put{$\cdot$} at 8.9 6.9

\put{$\vdots$} at 5.3 2.1
\put{$\vdots$} at 6.7 2.1

\put{$v$} at 6.5 5
\put{$\alpha_{1}$} at 1.5 8.5
\put{$\alpha_{2}$} at 2.5 9.5
\put{$\beta_{1}$} at 9.5 9.5
\put{$\beta_{2}$} at 10.5 8.5
\put{$\gamma_{1}$} at 6.7 0.3
\put{$\gamma_{2}$} at 5.3 0.3

\plot 1.7 7.8 3.7 5.8 /
\plot 1.7 7.8 1.8 7.9 / 
\plot 3.7 5.8 3.8 5.9 /

\plot 3.3 9.3 5.3 7.3 /
\plot 3.3 9.3 3.2 9.2 / 
\plot 5.3 7.3 5.2 7.2 /

\plot 6.7 7.3 8.7 9.3 /
\plot 6.7 7.3 6.8 7.2 / 
\plot 8.7 9.3 8.8 9.2 /

\plot 8.3 5.7 10.3 7.7 /
\plot 8.3 5.7 8.2 5.8 / 
\plot 10.3 7.7 10.2 7.8 /

\putrule from 4.9 3.3 to 4.9 .7
\putrule from 4.9 3.3 to 5.1 3.3
\putrule from 4.9 .7 to 5.1 .7

\putrule from 7.1 3.3 to 7.1 .7
\putrule from 7.1 3.3 to 6.9 3.3
\putrule from 7.1 .7 to 6.9 .7

\put{$l$} at 2.5 6.5
\put{$l$} at 4.5 8.5
\put{$l$} at 7.5 8.5
\put{$l$} at 9.5 6.5
\put{$l$} at 4.5 2
\put{$l$} at 7.5 2

\put{$W$} at 6 9.5
\put{Figure $8.$} at 6  -0.9
\endpicture 

\noindent
Note that $d(W) = 2\ell+2$.

\begin{theorem}
\label{q(W)}
Let $W$ be the $(3,\ell)$-whirl with $\ell \geq 2$.
Then
$$
q(W) \geq \displaystyle{\frac{9d(W)}{8}+ \frac{1}{2}} = d(W)+1+ \displaystyle\frac{\ell-1}{4} > d(W) +1.
$$
\end{theorem}

\noindent
{\bf Proof.}
Applying Proposition~\ref{path cover number} to the four paths illustrated in Figure 9, we conclude $p(W) =4$.

\beginpicture
\small
\setcoordinatesystem units <.700000cm, .700000cm>
\setplotarea x from -4.2 to 12, y from -1.5 to 10.5
\setlinear
\linethickness=0.50pt
\setplotsymbol ({\fiverm .})

\put {$\displaystyle\bullet$} at 2 8
\put {$\displaystyle\bullet$} at 2.5 7.5
\put {$\displaystyle\bullet$} at 3.5 6.5
\put {$\displaystyle\bullet$} at 4 6

\put {$\displaystyle\bullet$} at 3 9
\put {$\displaystyle\bullet$} at 3.5 8.5
\put {$\displaystyle\bullet$} at 4.5 7.5
\put {$\displaystyle\bullet$} at 5 7

\put {$\displaystyle\bullet$} at 5 6
\put {$\displaystyle\bullet$} at 6 4

\put {$\displaystyle\bullet$} at 9 9
\put {$\displaystyle\bullet$} at 8.5 8.5
\put {$\displaystyle\bullet$} at 7.5 7.5
\put {$\displaystyle\bullet$} at 7 7

\put {$\displaystyle\bullet$} at 10 8
\put {$\displaystyle\bullet$} at 9.5 7.5
\put {$\displaystyle\bullet$} at 8.5 6.5
\put {$\displaystyle\bullet$} at 8 6

\put {$\displaystyle\bullet$} at 7 6

\put {$\displaystyle\bullet$} at 6.7 3.3
\put {$\displaystyle\bullet$} at 6.7 2.5
\put {$\displaystyle\bullet$} at 6.7 1.5
\put {$\displaystyle\bullet$} at 6.7 .7

\put {$\displaystyle\bullet$} at 5.3 3.3
\put {$\displaystyle\bullet$} at 5.3 2.5
\put {$\displaystyle\bullet$} at 5.3 1.5
\put {$\displaystyle\bullet$} at 5.3 .7

\put {$\displaystyle\bullet$} at 6 5

\plot 2 8 2.5 7.5 /
\plot 3.5 6.5 4 6 /
\plot 4 6  5 6 /

\plot 3 9 3.5 8.5 /
\plot 4.5 7.5 5 7 /
\plot 5 7 5 6 /

\plot 9 9 8.5 8.5 /
\plot 7.5 7.5 7 7 /
\plot 7 7 7 6 /

\plot 10 8 9.5 7.5 /
\plot 8.5 6.5 8 6 /
\plot 8 6 7 6 /


\plot 6 4 5.3 3.3 /
\plot 5.3 3.3 5.3 2.5 /
\plot 5.3 1.5 5.3 .7 /

\plot 6 4 6.7 3.3 /
\plot 6.7 3.3 6.7 2.5 /
\plot 6.7 1.5 6.7 .7 /

\put{$\cdot$} at 3 7
\put{$\cdot$} at 2.9 7.1
\put{$\cdot$} at 3.1 6.9

\put{$\cdot$} at 4 8
\put{$\cdot$} at 3.9 8.1
\put{$\cdot$} at 4.1 7.9

\put{$\cdot$} at 8 8
\put{$\cdot$} at 8.1 8.1
\put{$\cdot$} at 7.9 7.9

\put{$\cdot$} at 9 7
\put{$\cdot$} at 9.1 7.1
\put{$\cdot$} at 8.9 6.9

\put{$\vdots$} at 5.3 2.1
\put{$\vdots$} at 6.7 2.1

\put{$v$} at 6.5 5
\put{$\alpha_{1}$} at 1.5 8.5
\put{$\alpha_{2}$} at 2.5 9.5
\put{$\beta_{1}$} at 9.5 9.5
\put{$\beta_{2}$} at 10.5 8.5
\put{$\gamma_{1}$} at 6.7 0.3
\put{$\gamma_{2}$} at 5.3 0.3

\put{$P_{1}$} at 2.5 6.5
\put{$P_{2}$} at 9.5 6.5
\put{$P_{3}$} at 7.5 2

\put{Figure $9.$} at 6  -0.9
\endpicture 

\noindent
Let $A \in S(W)$, and let $n_{j}$ be the number of the eigenvalues of $A$ with multiplicity $j$.
Then, since $p(W)=4$, Theorem~\ref{M=p} implies $n_{j} = 0$ for $j > 4$.
The three disjoint paths $P_{1}, P_{2}$ and $P_{3}$ in $W$ in Figure 9 cover all of the vertices of $W$ except $v$.
Thus, by Corollary~\ref{main cor},
\begin{equation}
\label{n4}
n_{4} \leq 1.
\end{equation}

Suppose that $\lambda \in \sigma(A)$ and $m_{A}(\lambda)=4$.
Consider the three disjoint paths $P_{1}, P_{2}$ and $P_{3}$ in $W$ illustrated in Figure 10.

\beginpicture
\small
\setcoordinatesystem units <.70000cm, .700000cm>
\setplotarea x from -4.2 to 12, y from -1.5 to 10.5
\setlinear
\linethickness=0.50pt
\setplotsymbol ({\fiverm .})

\put {$\displaystyle\bullet$} at 2 8
\put {$\displaystyle\bullet$} at 2.5 7.5
\put {$\displaystyle\bullet$} at 3.5 6.5
\put {$\displaystyle\bullet$} at 4 6

\put {$\displaystyle\bullet$} at 3 9
\put {$\displaystyle\bullet$} at 3.5 8.5
\put {$\displaystyle\bullet$} at 4.5 7.5
\put {$\displaystyle\bullet$} at 5 7

\put {$\displaystyle\bullet$} at 5 6
\put {$\displaystyle\bullet$} at 6 4

\put {$\displaystyle\bullet$} at 9 9
\put {$\displaystyle\bullet$} at 8.5 8.5
\put {$\displaystyle\bullet$} at 7.5 7.5
\put {$\displaystyle\bullet$} at 7 7

\put {$\displaystyle\bullet$} at 10 8
\put {$\displaystyle\bullet$} at 9.5 7.5
\put {$\displaystyle\bullet$} at 8.5 6.5
\put {$\displaystyle\bullet$} at 8 6

\put {$\displaystyle\bullet$} at 7 6

\put {$\displaystyle\bullet$} at 6.7 3.3
\put {$\displaystyle\bullet$} at 6.7 2.5
\put {$\displaystyle\bullet$} at 6.7 1.5
\put {$\displaystyle\bullet$} at 6.7 .7

\put {$\displaystyle\bullet$} at 5.3 3.3
\put {$\displaystyle\bullet$} at 5.3 2.5
\put {$\displaystyle\bullet$} at 5.3 1.5
\put {$\displaystyle\bullet$} at 5.3 .7

\put {$\displaystyle\bullet$} at 6 5

\plot 2 8 2.5 7.5 /
\plot 3.5 6.5 4 6 /

\plot 3 9 3.5 8.5 /
\plot 4.5 7.5 5 7 /
\plot 5 7 5 6 /

\plot 9 9 8.5 8.5 /
\plot 7.5 7.5 7 7 /
\plot 7 7 7 6 /

\plot 10 8 9.5 7.5 /
\plot 8.5 6.5 8 6 /
\plot 8 6 7 6 /

\plot 5 6  6 5 /
\plot 6 5 6 4 /

\plot 6 4 5.3 3.3 /
\plot 5.3 3.3 5.3 2.5 /
\plot 5.3 1.5 5.3 .7 /

\plot 6.7 3.3 6.7 2.5 /
\plot 6.7 1.5  6.7 .7 /

\put{$\cdot$} at 3 7
\put{$\cdot$} at 2.9 7.1
\put{$\cdot$} at 3.1 6.9

\put{$\cdot$} at 4 8
\put{$\cdot$} at 3.9 8.1
\put{$\cdot$} at 4.1 7.9

\put{$\cdot$} at 8 8
\put{$\cdot$} at 8.1 8.1
\put{$\cdot$} at 7.9 7.9

\put{$\cdot$} at 9 7
\put{$\cdot$} at 9.1 7.1
\put{$\cdot$} at 8.9 6.9

\put{$\vdots$} at 5.3 2.1
\put{$\vdots$} at 6.7 2.1

\put{$v$} at 6.5 5
\put{$\alpha_{1}$} at 1.5 8.5
\put{$\alpha_{2}$} at 2.5 9.5
\put{$\beta_{1}$} at 9.5 9.5
\put{$\beta_{2}$} at 10.5 8.5
\put{$\gamma_{2}$} at 5.3 0.3

\put{$P_{1}$} at 4.5 8.5
\put{$P_{2}$} at 9.5 6.5
\put{$P_{3}=\gamma_{1}$} at 8 2

\put{Figure $10.$} at 6  -0.9
\endpicture

\noindent
Then, by Corollary~\ref{main}, $\lambda$ is an eigenvalue of $A[W \setminus (P_{1} \cup P_{2} \cup P_{3})] = A[\alpha_{1}]$.
Similarly, $\lambda$ is an eigenvalue of $A[\alpha_{2}], A[\beta_{i}]$ and $A[\gamma_{i}]$ for $i=1,2$.
Since $\alpha_{i}, \beta_{i}$ and $\gamma_{i}$ are paths, Theorem~\ref{M=p} implies that
$\lambda$ is a simple eigenvalue of $A[\alpha_{i}], A[\beta_{i}]$ and $A[\gamma_{i}]$ for each $i=1,2$.

Next, suppose $\mu \in \sigma(A)$ and $m_{A}(\mu)=3$.
Let $P_{1}$ and  $P_{2}$ be the disjoint paths of $W$ illustrated in Figure 11.

\beginpicture
\small
\setcoordinatesystem units <.700000cm, .700000cm>
\setplotarea x from -4.2 to 12, y from -1.5 to 10.5
\setlinear
\linethickness=0.50pt
\setplotsymbol ({\fiverm .})

\put {$\displaystyle\bullet$} at 2 8
\put {$\displaystyle\bullet$} at 2.5 7.5
\put {$\displaystyle\bullet$} at 3.5 6.5
\put {$\displaystyle\bullet$} at 4 6

\put {$\displaystyle\bullet$} at 3 9
\put {$\displaystyle\bullet$} at 3.5 8.5
\put {$\displaystyle\bullet$} at 4.5 7.5
\put {$\displaystyle\bullet$} at 5 7

\put {$\displaystyle\bullet$} at 5 6
\put {$\displaystyle\bullet$} at 6 4

\put {$\displaystyle\bullet$} at 9 9
\put {$\displaystyle\bullet$} at 8.5 8.5
\put {$\displaystyle\bullet$} at 7.5 7.5
\put {$\displaystyle\bullet$} at 7 7

\put {$\displaystyle\bullet$} at 10 8
\put {$\displaystyle\bullet$} at 9.5 7.5
\put {$\displaystyle\bullet$} at 8.5 6.5
\put {$\displaystyle\bullet$} at 8 6

\put {$\displaystyle\bullet$} at 7 6

\put {$\displaystyle\bullet$} at 6.7 3.3
\put {$\displaystyle\bullet$} at 6.7 2.5
\put {$\displaystyle\bullet$} at 6.7 1.5
\put {$\displaystyle\bullet$} at 6.7 .7

\put {$\displaystyle\bullet$} at 5.3 3.3
\put {$\displaystyle\bullet$} at 5.3 2.5
\put {$\displaystyle\bullet$} at 5.3 1.5
\put {$\displaystyle\bullet$} at 5.3 .7

\put {$\displaystyle\bullet$} at 6 5

\plot 2 8 2.5 7.5 /
\plot 3.5 6.5 4 6 /

\plot 3 9 3.5 8.5 /
\plot 4.5 7.5 5 7 /
\plot 5 7 5 6 /

\plot 9 9 8.5 8.5 /
\plot 7.5 7.5 7 7 /
\plot 7 7 7 6 /

\plot 10 8 9.5 7.5 /
\plot 8.5 6.5 8 6 /
\plot 8 6 7 6 /

\plot 5 6  6 5 /
\plot 6 5 6 4 /

\plot 5.3 3.3 5.3 2.5 /
\plot 5.3 1.5 5.3 .7 /

\plot 6 4 6.7 3.3 /
\plot 6.7 3.3 6.7 2.5 /
\plot 6.7 1.5 6.7 .7 /

\put{$\cdot$} at 3 7
\put{$\cdot$} at 2.9 7.1
\put{$\cdot$} at 3.1 6.9

\put{$\cdot$} at 4 8
\put{$\cdot$} at 3.9 8.1
\put{$\cdot$} at 4.1 7.9

\put{$\cdot$} at 8 8
\put{$\cdot$} at 8.1 8.1
\put{$\cdot$} at 7.9 7.9

\put{$\cdot$} at 9 7
\put{$\cdot$} at 9.1 7.1
\put{$\cdot$} at 8.9 6.9

\put{$\vdots$} at 5.3 2.1
\put{$\vdots$} at 6.7 2.1

\put{$v$} at 6.5 5
\put{$\alpha_{1}$} at 1.5 8.5
\put{$\alpha_{2}$} at 2.5 9.5
\put{$\beta_{1}$} at 9.5 9.5
\put{$\beta_{2}$} at 10.5 8.5
\put{$\gamma_{1}$} at 6.7 0.3
\put{$\gamma_{2}$} at 5.3 0.3

\put{$P_{1}$} at 4.5 8.5
\put{$P_{2}$} at 9.5 6.5

\put{Figure $11.$} at 6  -0.9
\endpicture

\noindent
Then, by Corollary~\ref{main}, $\mu$ is an eigenvalue of $A[W \setminus (P_{1} \cup P_{2})]$ which is, up to permutation similarity,
equal to $A[\alpha_{1}] \oplus A[\gamma_{2}]$.
Thus, by Theorem~\ref{M=p}, $\mu$ is a simple eigenvalue of at least one of $A[\alpha_{1}]$, $A[\gamma_{2}]$.
Similarly, $\mu$ is a simple eigenvalue of $A[\alpha_{i}]$ or $A[\beta_{j}]$, $A[\alpha_{i}]$ or $A[\gamma_{j}]$,
and $A[\beta_{i}]$ or $A[\gamma_{j}]$ for $i, j \in \{1, 2\}$.
If $\mu$ is a simple eigenvalue of at most three matrices in
\begin{equation}
\label{display}\displaystyle\{ A[\alpha_{1}], A[\beta_{1}], A[\gamma_{1}], A[\alpha_{2}], A[\beta_{2}], A[\gamma_{2} ]\},
\end{equation}
then there is a pair of $A[\delta_{i}], A[\tau_{j}]$ for $\delta, \tau \in \{\alpha, \beta, \gamma\}$ and $\delta \ne \tau$
such that neither  $A[\delta_{i}]$ nor $ A[\tau_{j}]$ has $\mu$ as an eigenvalue.
Thus, $\mu$ is a simple eigenvalue of at least four matrices in $(\ref{display})$.

Let $A' = A[\alpha_{1}] \oplus A[\alpha_{2}] \oplus A[\beta_{1}] \oplus A[\beta_{2}] \oplus A[\gamma_{1}] \oplus A[\gamma_{2}]$.
We have shown (a) if $\lambda \in \sigma(A)$ and $m_{A}(\lambda)=4$, then $\lambda \in \sigma(A')$ and $m_{A'}(\lambda)=6$; and
(b) if $\mu \in \sigma(A)$ and $m_{A}(\mu)=3$, then $\mu \in \sigma(A')$ and $m_{A'}(\mu) \geq 4$.
Since the number of eigenvalues of a square matrix cannot exceed its order, by (a) and (b),
\begin{center}
$4n_{3}+6n_{4} \leq 6\ell$.
\end{center}

\noindent
Hence,
\begin{equation}
\label{n3}
n_{3} \leq \displaystyle{\frac{3\ell-3n_{4}}{2}}.
\end{equation}

By definition of $n_{j}$ and the fact that no eigenvalue of $A$ has multiplicity more than 4,
\begin{equation}
\label{n}
n=n_{1}+2n_{2}+3n_{3}+4n_{4} = n_{1}+2(n_{2}+n_{3}+n_{4})+n_{3}+2n_{4}.
\end{equation}
It is shown in \cite[Corollary 7]{JDS1} that the largest and smallest eigenvalues of $A$ are simple and hence, $n_{1} \geq 2$.
Since $n_{2}+n_{3}+n_{4}=q(A) - n_{1}$ and $n_{1} \geq 2$, by (\ref{n}),
\begin{equation}
\label{n=6l+4}
\begin{tabular}{rcl}
$n=6\ell+4$&$=$&$n_{1}+2(n_{2}+n_{3}+n_{4})+n_{3}+2n_{4}$ \\
           &$=$&$ n_{1}+2(q(A)-n_{1}) +n_{3} + 2 n_{4}$ \\
           &$=$& $2q(A) -n_{1} +n_{3} + 2 n_{4}$ \\
           &$\leq$& $2q(A) -2 +n_{3} + 2 n_{4}$.
\end{tabular}
\end{equation}

By (\ref{n4}), (\ref{n3}) and (\ref{n=6l+4}),  we have
\begin{equation}
\label{6l+4}
6\ell+4 \leq  2q(A)-2+ \frac{3\ell+n_{4}}{2} \leq 2q(A)-2+ \frac{3\ell+1}{2}.
\end{equation}

\noindent
By solving (\ref{6l+4}) for $q(A)$, and using $d(W)=2\ell+2$, we obtain
$$
q(A) \geq \frac{9\ell}{4}+ \frac{11}{4}= d(W)+1+ \frac{\ell-1}{4}= \frac{9d(W)}{8}+ \frac{1}{2}.
$$

\noindent
Since $A$ is an arbitrary matrix in $S(W)$, the result follows.
\hfill{\rule{2mm}{2mm}}
\vskip .1in

Next, we generalize Theorem~\ref{q(W)} to  whirls with more legs.
The \emph{$(k,\ell)$-whirl} $W$ for $k \geq 2$ and  $\ell \geq 1$ is the tree on $n=2k\ell+k+1$ vertices
with the axis vertex $v$ of degree $k$ and $2k$ legs $\{\alpha^{i}_{1}, \alpha^{i}_{2}\}_{i=1}^{k}$, each with $\ell$ vertices as illustrated in Figure 12.
Note that $d(W)=2\ell+2$.
If the number of vertices of each leg of $W$ is not specified, we say that $W$ is a \emph{$k$-whirl}.
If the numbers of the vertices of $2k$ legs are not necessarily equal, then $W$ is a \emph{generalized $k$-whirl}.

\beginpicture
\small
\setcoordinatesystem units <.700000cm, .700000cm>
\setplotarea x from -4.2 to 12, y from -1.5 to 10.5
\setlinear
\linethickness=0.50pt
\setplotsymbol ({\fiverm .})

\put {$\displaystyle\bullet$} at 2 8
\put {$\displaystyle\bullet$} at 2.5 7.5
\put {$\displaystyle\bullet$} at 3.5 6.5
\put {$\displaystyle\bullet$} at 4 6

\put {$\displaystyle\bullet$} at 3 9
\put {$\displaystyle\bullet$} at 3.5 8.5
\put {$\displaystyle\bullet$} at 4.5 7.5
\put {$\displaystyle\bullet$} at 5 7

\put {$\displaystyle\bullet$} at 5 6

\put {$\displaystyle\bullet$} at 9 9
\put {$\displaystyle\bullet$} at 8.5 8.5
\put {$\displaystyle\bullet$} at 7.5 7.5
\put {$\displaystyle\bullet$} at 7 7

\put {$\displaystyle\bullet$} at 10 8
\put {$\displaystyle\bullet$} at 9.5 7.5
\put {$\displaystyle\bullet$} at 8.5 6.5
\put {$\displaystyle\bullet$} at 8 6

\put {$\displaystyle\bullet$} at 7 6

\put {$\displaystyle\bullet$} at 4 4
\put {$\displaystyle\bullet$} at 3.5 3.5
\put {$\displaystyle\bullet$} at 2.5 2.5
\put {$\displaystyle\bullet$} at 2 2

\put {$\displaystyle\bullet$} at 5 3
\put {$\displaystyle\bullet$} at 4.5 2.5
\put {$\displaystyle\bullet$} at 3.5 1.5
\put {$\displaystyle\bullet$} at 3 1

\put {$\displaystyle\bullet$} at 5 4
\put {$\displaystyle\bullet$} at 6 5

\put {$\scriptscriptstyle\bullet$} at 7 3
\put {$\scriptscriptstyle\bullet$} at 8 3.5
\put {$\scriptscriptstyle\bullet$} at 8.5 4.5

\plot 2 8 2.5 7.5 /
\plot 3.5 6.5 4 6 /
\plot 4 6 5 6 /

\plot 3 9 3.5 8.5 /
\plot 4.5 7.5 5 7 /
\plot 5 7 5 6 /

\plot 9 9 8.5 8.5 /
\plot 7.5 7.5 7 7 /
\plot 7 7 7 6 /

\plot 10 8 9.5 7.5 /
\plot 8.5 6.5 8 6 /
\plot 8 6 7 6 /

\plot 5 6  6 5 /
\plot 6 5 7 6 /
\plot 6 5 5 4 /

\plot 4 4 5 4 /
\plot 4 4 3.5 3.5 /
\plot 2.5 2.5 2 2 /

\plot 5 4 5 3 /
\plot 5 3 4.5 2.5 /
\plot 3.5 1.5 3 1 /

\put{$\cdot$} at 3 7
\put{$\cdot$} at 2.9 7.1
\put{$\cdot$} at 3.1 6.9

\put{$\cdot$} at 4 8
\put{$\cdot$} at 3.9 8.1
\put{$\cdot$} at 4.1 7.9

\put{$\cdot$} at 8 8
\put{$\cdot$} at 8.1 8.1
\put{$\cdot$} at 7.9 7.9

\put{$\cdot$} at 9 7
\put{$\cdot$} at 9.1 7.1
\put{$\cdot$} at 8.9 6.9

\put{$\cdot$} at 3 3
\put{$\cdot$} at 3.1 3.1
\put{$\cdot$} at 2.9 2.9

\put{$\cdot$} at 4 2
\put{$\cdot$} at 4.1 2.1
\put{$\cdot$} at 3.9 1.9

\put{$v$} at 6.5 5
\put{$v_{1}$} at 5.5 6
\put{$v_{2}$} at 6.5 6
\put{$v_{k}$} at 5.5 4

\put{$\alpha^{1}_{1}$} at 1.5 8.5
\put{$\alpha^{1}_{2}$} at 2.5 9.5
\put{$\alpha^{2}_{1}$} at 9.5 9.5
\put{$\alpha^{2}_{2}$} at 10.5 8.5
\put{$\alpha^{k}_{2}$} at 1.5 1.5
\put{$\alpha^{k}_{1}$} at 2.5 .5

\put{$W$} at 6 9.5
\put{Figure $12.$} at 6  -0.9
\endpicture 

In order to prove a generalization of the result in Theorem~\ref{q(W)} for the $(k,\ell)$-whirl,
we first prove the following.

\begin{lemma}
\label{n's}
Suppose that $W$ is a generalized $k$-whirl $(k \geq 2)$ on $n$ vertices with $2k$ legs $\{\alpha^{i}_{1}, \alpha^{i}_{2}\}_{i=1}^{k}$
such that each leg has at least one vertex.
Let $A \in S(W)$, and let $A'$ be the direct sum of $A[\alpha^{i}_{j}]$ for all $i \in \{1, \ldots, k\}$, $j \in \{1,2\}$.
If $n_{r}$ denotes the number of eigenvalues of $A$ with multiplicity $r$, then the following hold:

\begin{itemize}
\item[\rm(a)] $n_{k+1} \leq 1$, and $n_{j}=0$ for $j \geq k+2$;
\item[\rm(b)] If $\lambda \in \sigma(A)$ and $m_{A}(\lambda) = k+1$, then
$\lambda$ is a simple eigenvalue of $A[\alpha^{i}_{j}]$ for all $i \in \{1, \ldots, k\}$, $j \in \{1,2\}$,
and $m_{A'}(\lambda) = 2k$;
\item[\rm(c)] If $\mu \in \sigma(A)$ and $m_{A}(\mu) =k$, then, for all $i \ne s$ and all $j$ and $t$,
$\mu$ is a simple eigenvalue of at least one of $A[\alpha^{i}_{j}]$, $A[\alpha^{s}_{t}]$, and $m_{A'}(\mu) \geq 2k-2$; and
\item[\rm(d)] $(2k-2)n_{k} + (2k)n_{k+1} \leq n-(k+1)$.
\end{itemize}
\end{lemma}

\noindent
{\bf Proof.}
{\bf (a)} By applying Proposition~\ref{path cover number} to the $k+1$ paths illustrated in Figure 13, we conclude $p(W) =k+1$.

\beginpicture
\small
\setcoordinatesystem units <.700000cm, .700000cm>
\setplotarea x from -4.2 to 12, y from -1.5 to 10.5
\setlinear
\linethickness=0.50pt
\setplotsymbol ({\fiverm .})

\put {$\displaystyle\bullet$} at 2 8
\put {$\displaystyle\bullet$} at 2.5 7.5
\put {$\displaystyle\bullet$} at 3.5 6.5
\put {$\displaystyle\bullet$} at 4 6

\put {$\displaystyle\bullet$} at 3 9
\put {$\displaystyle\bullet$} at 3.5 8.5
\put {$\displaystyle\bullet$} at 4.5 7.5
\put {$\displaystyle\bullet$} at 5 7

\put {$\displaystyle\bullet$} at 5 6

\put {$\displaystyle\bullet$} at 9 9
\put {$\displaystyle\bullet$} at 8.5 8.5
\put {$\displaystyle\bullet$} at 7.5 7.5
\put {$\displaystyle\bullet$} at 7 7

\put {$\displaystyle\bullet$} at 10 8
\put {$\displaystyle\bullet$} at 9.5 7.5
\put {$\displaystyle\bullet$} at 8.5 6.5
\put {$\displaystyle\bullet$} at 8 6

\put {$\displaystyle\bullet$} at 7 6

\put {$\displaystyle\bullet$} at 4 4
\put {$\displaystyle\bullet$} at 3.5 3.5
\put {$\displaystyle\bullet$} at 2.5 2.5
\put {$\displaystyle\bullet$} at 2 2

\put {$\displaystyle\bullet$} at 5 3
\put {$\displaystyle\bullet$} at 4.5 2.5
\put {$\displaystyle\bullet$} at 3.5 1.5
\put {$\displaystyle\bullet$} at 3 1

\put {$\displaystyle\bullet$} at 5 4
\put {$\displaystyle\bullet$} at 6 5

\put {$\scriptscriptstyle\bullet$} at 7 3
\put {$\scriptscriptstyle\bullet$} at 8 3.5
\put {$\scriptscriptstyle\bullet$} at 8.5 4.5

\plot 2 8 2.5 7.5 /
\plot 3.5 6.5 4 6 /
\plot 4 6 5 6 /

\plot 3 9 3.5 8.5 /
\plot 4.5 7.5 5 7 /
\plot 5 7 5 6 /

\plot 9 9 8.5 8.5 /
\plot 7.5 7.5 7 7 /
\plot 7 7 7 6 /

\plot 10 8 9.5 7.5 /
\plot 8.5 6.5 8 6 /
\plot 8 6 7 6 /


\plot 4 4 5 4 /
\plot 4 4 3.5 3.5 /
\plot 2.5 2.5 2 2 /

\plot 5 4 5 3 /
\plot 5 3 4.5 2.5 /
\plot 3.5 1.5 3 1 /

\put{$\cdot$} at 3 7
\put{$\cdot$} at 2.9 7.1
\put{$\cdot$} at 3.1 6.9

\put{$\cdot$} at 4 8
\put{$\cdot$} at 3.9 8.1
\put{$\cdot$} at 4.1 7.9

\put{$\cdot$} at 8 8
\put{$\cdot$} at 8.1 8.1
\put{$\cdot$} at 7.9 7.9

\put{$\cdot$} at 9 7
\put{$\cdot$} at 9.1 7.1
\put{$\cdot$} at 8.9 6.9

\put{$\cdot$} at 3 3
\put{$\cdot$} at 3.1 3.1
\put{$\cdot$} at 2.9 2.9

\put{$\cdot$} at 4 2
\put{$\cdot$} at 4.1 2.1
\put{$\cdot$} at 3.9 1.9

\put{$v$} at 6.5 5
\put{$v_{1}$} at 5.5 6
\put{$v_{2}$} at 6.5 6
\put{$v_{k}$} at 5.5 4

\put{$\alpha^{1}_{1}$} at 1.5 8.5
\put{$\alpha^{1}_{2}$} at 2.5 9.5
\put{$\alpha^{2}_{1}$} at 9.5 9.5
\put{$\alpha^{2}_{2}$} at 10.5 8.5
\put{$\alpha^{k}_{2}$} at 1.5 1.5
\put{$\alpha^{k}_{1}$} at 2.5 .5

\put{$P_{1}$} at 2.5 6.5 
\put{$P_{2}$} at 9.5 6.5
\put{$P_{k}$} at 2.5 3.5

\put{Figure $13.$} at 6  -0.9
\endpicture 

\noindent
Hence, Theorem~\ref{M=p} implies that there is no eigenvalue of $A$ with multiplicity $k+2$ or more, that is,
$n_{j} =0$ for $j \geq k+2$.
Furthermore, since the $k$ disjoint paths $P_{1}, \ldots, P_{k}$ in $W$ in Figure 13 cover all of the vertices of $W$ except $v$,
Corollary~\ref{main cor} implies that there exists at most one eigenvalue of $A$ with multiplicity $k+1$, that is, $n_{k+1} \leq 1$.
This proves (a).

\medskip
\noindent
{\bf (b)} Suppose $\lambda \in \sigma(A)$ and $m_{A}(\lambda) = k+1$.
Consider the following $k$ disjoint paths in $W$:
$P_{1} = \alpha^{1}_{2}(v_{1}$---$v$---$v_{k}) \alpha^{k}_{2}$, $P_{i} = \alpha^{i}_{1} v_{i} \alpha^{i}_{2}$ for $i=2, \ldots, k-1$, and
$P_{k} = \alpha^{k}_{1}$.
Then $W \setminus (P_{1} \cup \cdots \cup P_{k})$ is $\alpha^{1}_{1}$ and hence, Corollary~\ref{main} implies $\lambda \in A[\alpha^{1}_{1}]$.
Since $\alpha^{1}_{1}$ is a path, by Theorem~\ref{M=p}, $\lambda$ is a simple eigenvalue of $A[\alpha^{1}_{1}]$.
Similarly, $\lambda$ is a simple eigenvalue of $A[\alpha^{1}_{2}]$ and $A[\alpha^{i}_{j}]$ for each $i=2, \ldots, k$
and $j=1,2$.
Therefore, $m_{A'}(\lambda) = 2k$, and (b) holds.

\medskip
\noindent
{\bf (c)} Suppose $\mu \in \sigma(A)$ and $m_{A}(\mu)=k$.
Let $P_{1} = \alpha^{1}_{2}(v_{1}$---$v$---$v_{k}) \alpha^{k}_{2}$, and $P_{i} = \alpha^{i}_{1} v_{i} \alpha^{i}_{2}$ for $i=2, \ldots, k-1$.
Then $A[W \setminus (P_{1} \cup \cdots \cup P_{k-1})]$ is, up to permutation similarity, equal to $A[\alpha^{1}_{1}] \oplus A[\alpha^{k}_{1}]$.
Thus, by Corollary~\ref{main}, $\mu$ is an eigenvalue of $A[\alpha^{1}_{1}] \oplus A[\alpha^{k}_{1}]$ and hence,
by Theorem~\ref{M=p}, $\mu$ is a simple eigenvalue of at least one of $A[\alpha^{1}_{1}]$, $A[\alpha^{k}_{1}]$.
Similarly, $\mu$ is a simple eigenvalue of at least one of $A[\alpha^{i}_{j}]$, $A[\alpha^{s}_{t}]$
for all $i, s \in \{1, \ldots, k\}$, $i \ne s$ and $j, t \in \{1,2\}$.
If $\mu$ is a simple eigenvalue of at most $2k-3$ of the $2k$ matrices in $\{A[\alpha^{i}_{1}], A[\alpha^{i}_{2}]\}_{i=1}^{k}$, then
there is a pair of $A[\alpha^{i}_{j}], A[\alpha^{s}_{t}]$ for $i \ne s$ such that none of
$A[\alpha^{i}_{j}], A[\alpha^{s}_{t}]$ has $\mu$ as an eigenvalue.
Thus, $\mu$ is a simple eigenvalue of at least $2k-2$ matrices in $\{A[\alpha^{i}_{1}], A[\alpha^{i}_{2}]\}_{i=1}^{k}$.
Therefore, $m_{A'}(\mu) \geq 2k-2$, and (c) holds.

\medskip
\noindent
{\bf (d)} The order of $A'$ is $n-(k+1)$.
Since the number of eigenvalues of $A'$ cannot exceed the order of $A'$,
$(2k-2)n_{k} + (2k)n_{k+1} \leq n-(k+1)$.
\hfill{\rule{2mm}{2mm}}

\begin{theorem}
\label{general q(W)}
Let $W$ be the $(k,\ell)$-whirl with $k \geq 3$ and  $ \ell \geq 2$.
Then
$$
q(W) \geq d(W)+1 + \displaystyle\frac{(k-2)(\ell-1)}{(k-1)^{2}} > d(W)+1.
$$
\end{theorem}

\noindent
{\bf Proof.}
Let $A \in S(W)$, and let $n_{j}$ be the number of eigenvalues of $A$ with multiplicity $j$.
By Lemma~\ref{n's} (a),
$$
n=2k\ell+k+1 = \displaystyle\sum_{i=1}^{k+1}{i \cdot n_{i}}.
$$
Moreover,
\begin{equation}
\label{2kl+k+1}
\begin{tabular}{rcl}
$2k\ell+k+1$&$ =$& $\displaystyle\sum_{i=1}^{k+1}{i \cdot n_{i}}$ \\
        &$=$&$ n_{1} + 2n_{2}+ \cdots + kn_{k} + (k+1)n_{k+1}$ \\
        &$\leq $& $2 + (k-1)(n_{1}-2 +n_{2}+\cdots + n_{k} + n_{k+1}) + n_{k} + 2n_{k+1}$ \\
        &$=$& $2 + (k-1)(q(A)-2) + n_{k} + 2n_{k+1}$.
\end{tabular}
\end{equation}

Next, since $n=2k\ell+k+1$, Lemma~\ref{n's} (d) implies
\begin{equation}
\label{leq 2kl}
(2k-2)n_{k} + (2k)n_{k+1} \leq 2k\ell.
\end{equation}
\noindent
Thus, by solving (\ref{leq 2kl}) for $n_{k}$, we have
$$
n_{k} \leq \displaystyle\frac{k\ell-kn_{k+1}}{k-1}.
$$
Furthermore, by (\ref{2kl+k+1}),
\begin{eqnarray*}
2k\ell+k+1 &\leq & 2+(k-1)(q(A)-2)+ \displaystyle\frac{k\ell-kn_{k+1}}{k-1} + 2n_{k+1} \\
        & =   & (k-1)q(A)-2k+4+ \displaystyle\frac{k\ell+(k-2)n_{k+1}}{k-1}.
\end{eqnarray*}
By Lemma~\ref{n's} (a), $n_{k+1} \leq 1$ and hence,
\begin{equation}
\label{2kl+3k-3}
2k\ell+3k-3 \leq (k-1)q(A) + \displaystyle\frac{k\ell+k-2}{k-1}.
\end{equation}
By solving (\ref{2kl+3k-3}) for $q(A)$, we have
\begin{eqnarray*}
q(A)  &\geq & \displaystyle\frac{[(2k\ell+3k-3)(k-1)]-(k\ell+k-2)}{(k-1)^{2}} \\
      & =   & \displaystyle\frac{2k^{2}\ell+3k^{2}-3k\ell-7k+5}{(k-1)^{2}} \\
      & =   & 2\ell+3 + \displaystyle\frac{(k-2)(\ell-1)}{(k-1)^{2}}.
\end{eqnarray*}
Since $d(W)=2\ell+2$, we have
\begin{equation}
\label{q(A)2}
q(A) \geq d(W)+1+ \displaystyle\frac{(k-2)(\ell-1)}{(k-1)^{2}}.
\end{equation}
Since $A$ is an arbitrary matrix in $S(W)$, (\ref{q(A)2}) implies
$$
q(W) \geq d(W)+1+ \displaystyle\frac{(k-2)(\ell-1)}{(k-1)^{2}} > d(W)+1
$$
for $k \geq 3$ and  $\ell \geq 2$.
\hfill{\rule{2mm}{2mm}}

\section{$q(G)$ for a Class of Connected Graphs $G$}
In this section we apply the Smith Normal Form approach to a class of connected graphs $G$, and
find a lower bound on $q(G)$.

Let $G$ be a connected graph on $n$ vertices with 6 legs $L_{i}$, each with $\ell$ vertices for $\ell \geq 1$, (see Figure 14) such that
$H$ is a connected graph on $m$ vertices containing $v_{1}, v_{2}, v_{3}$, and
for all $s$, $t$, and $k$ with $\{s,t,k\}=\{1,2,3\}$
there exists a unique shortest path from $v_{s}$ to $v_{t}$ which does not pass $v_{k}$.
We use $i_{u}, j_{v}$ to denote the pendant vertices of the 6 legs (see Figure 14).
   
\beginpicture
\small
\setcoordinatesystem units <.700000cm, .700000cm>
\setplotarea x from -3.5 to 14, y from -1 to 10.5
\setlinear
\linethickness=0.50pt
\setplotsymbol ({\fiverm .})

\put {$\displaystyle\bullet$} at 1 7.5
\put {$\displaystyle\bullet$} at 2 7.5
\put {$\displaystyle\bullet$} at 4 7.5
\put {$\displaystyle\bullet$} at 5 7.5

\put {$\displaystyle\bullet$} at 6 8

\put {$\displaystyle\bullet$} at 1 8.5
\put {$\displaystyle\bullet$} at 2 8.5
\put {$\displaystyle\bullet$} at 4 8.5
\put {$\displaystyle\bullet$} at 5 8.5

\put {$\displaystyle\bullet$} at 13 8.5
\put {$\displaystyle\bullet$} at 12 8.5
\put {$\displaystyle\bullet$} at 10 8.5
\put {$\displaystyle\bullet$} at 9 8.5

\put {$\displaystyle\bullet$} at 8 8

\put {$\displaystyle\bullet$} at 13 7.5
\put {$\displaystyle\bullet$} at 12 7.5
\put {$\displaystyle\bullet$} at 10 7.5
\put {$\displaystyle\bullet$} at 9 7.5

\put {$\displaystyle\bullet$} at 7.5 2
\put {$\displaystyle\bullet$} at 7.5 3
\put {$\displaystyle\bullet$} at 7.5 5
\put {$\displaystyle\bullet$} at 7.5 6

\put {$\displaystyle\bullet$} at 7 7

\put {$\displaystyle\bullet$} at 6.5 2
\put {$\displaystyle\bullet$} at 6.5 3
\put {$\displaystyle\bullet$} at 6.5 5
\put {$\displaystyle\bullet$} at 6.5 6

\put {$\cdots$} at 3 7.5
\put {$\cdots$} at 3 8.5

\put {$\cdots$} at 11 8.5
\put {$\cdots$} at 11 7.5

\put {$\vdots$} at 7.5 4
\put {$\vdots$} at 6.5 4

\circulararc 360 degrees from 7 9 center at 7 8

\plot 1 7.5 2 7.5 /
\plot 4 7.5 5 7.5 /
\plot 5 7.5 6 8 /

\plot 1 8.5 2 8.5 /
\plot 4 8.5 5 8.5 /
\plot 5 8.5 6 8 /

\plot 13 8.5 12 8.5 /
\plot 10 8.5 9 8.5 /
\plot 9 8.5 8 8 /

\plot 13 7.5 12 7.5 /
\plot 10 7.5 9 7.5 /
\plot 9 7.5 8 8 /

\plot 7.5 2 7.5 3 /
\plot 7.5 5 7.5 6 /
\plot 7.5 6 7 7 /

\plot 6.5 2 6.5 3 /
\plot 6.5 5 6.5 6 /
\plot 6.5 6 7 7 /

\put {$L_{1}$} at 3 6.9
\put {$L_{2}$} at 3 9.1

\put {$L_{3}$} at 11 9.1
\put {$L_{4}$} at 11 6.9

\put {$L_{5}$} at 8.1 4
\put {$L_{6}$} at 5.9 4

\put{$i_{1}$} at .5 7.5
\put{$j_{1}$} at .5 8.5

\put{$i_{2}$} at 13.5 8.5
\put{$j_{2}$} at 13.5 7.5

\put{$i_{3}$} at 7.5 1.5
\put{$j_{3}$} at 6.5 1.5

\put{$v_{1}$} at 5.8 8.5
\put{$v_{2}$} at 8.2 8.5
\put{$v_{3}$} at 7.5 6.7

\put{$H$} at 7 8

\put{Figure $14.$} at 7 0.1
\endpicture

\noindent
Using the Smith Normal Form approach, we show the following.

\begin{theorem}
\label{q(G)}
Let $G$ be the connected graph described in Figure $14$. Then
$$
q(G) \geq \frac{9\ell}{4} -2m + \frac{15}{2}.
$$
\end{theorem}

\noindent
{\bf Proof.}
Let $A \in S(G)$, and $\lambda_{1}, \ldots, \lambda_{q}$ be the distinct eigenvalues of $A$, and
let $n_{j}$ be the number of eigenvalues of $A$ with multiplicity $j$.
We consider $xI-A$.
Note that $D(xI-A)$ has a loop at each vertex.

Now we compute the determinant of an $(n-3)$ by $(n-3)$ submatrix of $xI-A$, and find
an upper bound on $\displaystyle\sum_{\lambda_{j}: m_{A}(\lambda_{j}) \geq 4}{(m_{A}(\lambda_{j}) -2)}$.
Let $M$ be the matrix obtained from $xI-A$ by replacing column $i_{s}$ by $e_{j_{s}}$ for each $s=1, 2,3$.
In terms of digraphs, $D(M)$ is obtained from $D(xI-A)$ by deleting all incoming arcs to vertices $i_{1}, i_{2}, i_{3}$, and
inserting the arcs $(j_{1}, i_{1}), (j_{2}, i_{2}), (j_{3}, i_{3})$, each with weight 1.
Hence, $(j_{s}, i_{s})$ is the unique arc of $D(M)$ ending at vertex $i_{s}$ for each $s=1, 2, 3$.
Moreover, if there exists a directed walk from $i_{s}$ to $j_{s}$ in $D(M)$ which is not
the directed path $P_{s}= (i_{s},\ldots, v_{s}, \ldots,  j_{s})$, then the directed walk has a loop or
it repeats $v_{s}$ at least twice.
Thus, $\beta_{s}= (P_{s}, (j_{s},i_{s}))$ is the unique directed cycle containing $(j_{s}, i_{s})$ for each $s=1, 2, 3$.
Therefore, each set of disjoint directed cycles of $D(M)$ that cover every vertex consists of the directed cycles $\beta_{s}$'s along with
 disjoint directed cycles covering every vertex of $H \setminus U$ where $U=\{v_{1}, v_{2}, v_{3}\}$.
By (\ref{digraph det}), this implies that
\begin{equation}
\label{det M}
\det M = \pm \displaystyle\prod_{s=1}^{3} \mbox{wt}(P_{s}) \cdot \det M[H \setminus U].
\end{equation}

\noindent
Note that $M[H \setminus U]= (xI-A)[H \setminus U]$, and $\mbox{wt}(P_{s})$ is a nonzero constant for each $s=1,2,3$.
By  Laplace expansion of the determinant along the columns $i_{1}, i_{2}, i_{3}$ of $M$,
$\det M  = \pm \det[(xI-A)(\{j_{1}, j_{2}, j_{3}\}, \{i_{1}, i_{2}, i_{3} \})]$.
Thus, by (\ref{det M}),
$$
\det[(xI-A)(\{j_{1}, j_{2}, j_{3} \}, \{ i_{1}, i_{2}, i_{3} \})] = c \cdot \det((xI-A)[H \setminus U])
$$
for some nonzero constant $c$.

Since $(xI-A)(\{j_{1}, j_{2}, j_{3}\}, \{ i_{1}, i_{2}, i_{3} \})$ is an $n-3$ by $n-3$ submatrix of $xI-A$,
we have
$$
\Delta_{n-3}(x) | \det((xI-A)[H \setminus U]).
$$
Thus, by Theorem~\ref{factor main},
if $\lambda$ is an eigenvalue of $A$ with $m_{A}(\lambda) \geq 4$, then
$\lambda$ is a zero of $(xI-A)[H \setminus U]$ with multiplicity $m_{A}(\lambda)-3$ or more.
Since the order of $(xI-A)[H \setminus U]$ is $m-3$, we have
\begin{equation}
\label{leq m-3}
\displaystyle\sum_{\lambda_{j}: m_{A}(\lambda_{j}) \geq 4}{1}~
\leq \displaystyle\sum_{\lambda_{j}: m_{A}(\lambda_{j}) \geq 4}{(m_{A}(\lambda_{j}) -3)}
\leq m-3.
\end{equation}
Thus, by (\ref{leq m-3}),
\begin{equation}
\label{leq 2m-6}
\displaystyle\sum_{\lambda_{j}: m_{A}(\lambda_{j}) \geq 4}{(m_{A}(\lambda_{j}) -2)}~
= \displaystyle\sum_{\lambda_{j}: m_{A}(\lambda_{j}) \geq 4}{(m_{A}(\lambda_{j}) -3)} ~+  \displaystyle\sum_{\lambda_{j}: m_{A}(\lambda_{j}) \geq 4}{1}
~ \leq ~ 2m-6.
\end{equation}

Next, we compute the determinants of $(n-2)$ by $(n-2)$ submatrices of $xI-A$ and thereby, find
an upper bound on $n_{3}$.
Let $N$ be the matrix obtained from $xI-A$ by replacing column $a_{s}$ by $e_{b_{t}}$, and
column $i_{k}$ by $e_{j_{k}}$ where $a, b \in \{i,j\}$ and $\{s, t, k\} = \{1,2,3\}$.
In terms of digraphs, $D(N)$ is obtained from $D(xI-A)$ by deleting all incoming arcs to $a_{s}, i_{k}$, and
inserting $(b_{t}, a_{s}), (j_{k}, i_{k})$, each with weight 1 (see Figure 15).

\beginpicture
\small
\setcoordinatesystem units <.700000cm, .700000cm>
\setplotarea x from -3.5 to 14, y from -1 to 10.5
\setlinear
\linethickness=0.50pt
\setplotsymbol ({\fiverm .})

\put {$\displaystyle\bullet$} at 1 7.5
\put {$\displaystyle\bullet$} at 2 7.5
\put {$\displaystyle\bullet$} at 4 7.5
\put {$\displaystyle\bullet$} at 5 7.5

\put {$\displaystyle\bullet$} at 6 8

\put {$\displaystyle\bullet$} at 1 8.5
\put {$\displaystyle\bullet$} at 2 8.5
\put {$\displaystyle\bullet$} at 4 8.5
\put {$\displaystyle\bullet$} at 5 8.5

\put {$\displaystyle\bullet$} at 13 8.5
\put {$\displaystyle\bullet$} at 12 8.5
\put {$\displaystyle\bullet$} at 10 8.5
\put {$\displaystyle\bullet$} at 9 8.5

\put {$\displaystyle\bullet$} at 8 8

\put {$\displaystyle\bullet$} at 13 7.5
\put {$\displaystyle\bullet$} at 12 7.5
\put {$\displaystyle\bullet$} at 10 7.5
\put {$\displaystyle\bullet$} at 9 7.5

\put {$\displaystyle\bullet$} at 7.5 2
\put {$\displaystyle\bullet$} at 7.5 3
\put {$\displaystyle\bullet$} at 7.5 5
\put {$\displaystyle\bullet$} at 7.5 6

\put {$\displaystyle\bullet$} at 7 7

\put {$\displaystyle\bullet$} at 6.5 2
\put {$\displaystyle\bullet$} at 6.5 3
\put {$\displaystyle\bullet$} at 6.5 5
\put {$\displaystyle\bullet$} at 6.5 6

\put {$\cdots$} at 3 7.5
\put {$\cdots$} at 3 8.5

\put {$\cdots$} at 11 8.5
\put {$\cdots$} at 11 7.5

\put {$\vdots$} at 7.5 4
\put {$\vdots$} at 6.5 4

\circulararc 360 degrees from 7 9 center at 7 8

\plot 1 7.5 2 7.5 /
\plot 4 7.5 5 7.5 /
\plot 5 7.5 6 8 /

\plot 1 8.5 2 8.5 /
\plot 4 8.5 5 8.5 /
\plot 5 8.5 6 8 /

\plot 13 8.5 12 8.5 /
\plot 10 8.5 9 8.5 /
\plot 9 8.5 8 8 /

\plot 13 7.5 12 7.5 /
\plot 10 7.5 9 7.5 /
\plot 9 7.5 8 8 /

\plot 7.5 2 7.5 3 /
\plot 7.5 5 7.5 6 /
\plot 7.5 6 7 7 /

\plot 6.5 2 6.5 3 /
\plot 6.5 5 6.5 6 /
\plot 6.5 6 7 7 /

\plot 13.25 8 13.15 7.9 /
\plot 13.25 8 13.35 7.9 /

\plot 12.5 8.5 12.6 8.6 /
\plot 12.5 8.5 12.6 8.4 /

\plot 1.5 7.5 1.4 7.6 /
\plot 1.5 7.5 1.4 7.4 /

\plot 2 4 2.1 4 /
\plot 2 4 2 3.9 /

\setquadratic
\plot 6.5 2 2 4 1 7.5 /
\plot 13 7.5 13.25 8 13 8.5 /

\put {$L$} at 3 9


\put {$R$} at 8 4

\put{$a_{s}$} at .5 7.5

\put{$i_{k}$} at 13.5 8.5
\put{$j_{k}$} at 13.5 7.5

\put{$b_{t}$} at 6.5 1.5

\put{$v_{s}$} at 5.8 8.5
\put{$v_{k}$} at 8.2 8.5
\put{$v_{t}$} at 7.5 6.7

\put{$H$} at 7 8

\put{Figure $15.$} at 7 .1
\endpicture

In Figure 15, edges represent directed 2-cycles.
Therefore, $D(N)$ has exactly one arc, namely $(b_{t}, a_{s})$, ending at $a_{s}$ and exactly one arc, namely $(j_{k}, i_{k})$, ending at $i_{k}$.
Moreover, each directed path from $a_{s}$ to $b_{t}$ passes through $v_{s}$ and $v_{t}$.
Let $Q^{k}_{1}, \ldots, Q^{k}_{d_{k}}$ be the directed paths from $v_{s}$ to $v_{t}$ not containing vertex $v_{k}$, and
$Q^{k}_{1}$ be the unique shortest directed path from $v_{s}$ to $v_{t}$ that does not go through $v_{k}$.
Then, for each $g=1, \ldots, d_{k}$, $(P_{a_{s} \to v_{s}}, Q^{k}_{g}, P_{v_{t} \to b_{t}})$ is
a directed path from $a_{s}$ to $b_{t}$ not containing vertex $v_{k}$, and
$(P_{a_{s} \to v_{s}}, Q^{k}_{1}, P_{v_{t} \to b_{t}})$ is the unique shortest directed path from $a_{s}$ to $b_{t}$ that does not go through vertex $v_{k}$.
Note that there exists exactly one directed cycle $\alpha$ of $D(N)$ containing $(j_{k}, i_{k})$,
$\alpha = (P_{i_{k} \to v_{k}}, P_{v_{k} \to j_{k}}, (j_{k}, i_{k}))$.
Thus, the directed cycles $\beta^{k}_{g}= (P_{a_{s} \to v_{s}}, Q^{k}_{g}, P_{v_{t} \to b_{t}}, (b_{t}, a_{s}))$ $(g=1, \ldots, d_{k})$ are the only directed cycles containing $(b_{t}, a_{s})$, which are disjoint from $\alpha$.
Let $H' = H \setminus \{v_{k} \}$.
Then, by (\ref{digraph det}),
\begin{equation}
\label{det N1}
\det N  = \mbox{swt}(\alpha) \cdot \det N[L] \cdot \det N[R] \cdot
\displaystyle\sum_{g=1}^{d_{k}}{\mbox{swt}(\beta^{k}_{g}) \cdot \det N[H' \setminus \beta^{k}_{g}]},
\end{equation}
where $N[H' \setminus \beta^{k}_{g}]$ is the principal submatrix of $N$ whose rows and columns correspond to $V(H') \setminus V(\beta^{k}_{g})$.

 Since $H' \setminus \beta^{k}_{1}$ is the unique one among $(H' \setminus \beta^{k}_{g})$'s having the largest number of
vertices and $\mbox{swt}(\beta^{k}_{g})$'s are nonzero constants,
$\displaystyle \sum_{g=1}^{d_{k}}\mbox{swt}(\beta^{k}_{g})\cdot \det N[H' \setminus \beta^{k}_{g}]$ is a nonzero polynomial
of degree $|V(H') \setminus V(\beta^{k}_{1})|$.

Note that $\mbox{swt}(\beta^{k}_{g}) = (-1)^{v(P_{a_{s} \to v_{s}}) + v(Q^{k}_{g}) + v(P_{v_{t} \to b_{t}}) -3}
\mbox{wt}(P_{a_{s} \to v_{s}}) \cdot \mbox{wt}(Q^{k}_{g}) \cdot \mbox{wt}(P_{v_{t} \to b_{t}})$, and
$V(H') \setminus V(\beta^{k}_{g}) = V(H') \setminus V(Q^{k}_{g})$, where
$v(P_{a_{s} \to v_{s}})$, $v(Q^{k}_{g})$ and $v(P_{v_{t} \to b_{t}})$ are the numbers of the vertices of those directed paths.
Hence, from (\ref{det N1}),

$$
\begin{tabular}{rcl}
$\det N$ &=& $\mbox{swt}(\alpha) \cdot \det N[L] \cdot \det N[R] \cdot (-1)^{v(P_{a_{s} \to v_{s}}) + v(P_{v_{t} \to b_{t}}) -3}
\mbox{wt}(P_{a_{s} \to v_{s}}) \cdot \mbox{wt}(P_{v_{t} \to b_{t}})$  \\
& &$\cdot \displaystyle\sum_{g=1}^{d_{k}}(-1)^{v(Q^{k}_{g})}\mbox{wt}(Q^{k}_{g}) \cdot \det N[H' \setminus Q^{k}_{g}]$.
\end{tabular}
$$

\noindent
Since $\mbox{swt}(\alpha)$, $\mbox{wt}(P_{a_{s} \to v_{s}})$ and $\mbox{wt}(P_{v_{t} \to b_{t}})$ are nonzero constants,
\begin{equation}
\label{det N2}
\det N = c \cdot \det N[L] \cdot \det N[R] \cdot \displaystyle\sum_{g=1}^{d_{k}}(-1)^{v(Q^{k}_{g})}\mbox{wt}(Q^{k}_{g}) \cdot \det N[H' \setminus Q^{k}_{g}],
\end{equation}
for some nonzero constant $c$.

Note that $\det N[L]= \det ((xI-A)[L])$, and $\det N[R]= \det((xI-A)[R])$.
By Laplace expansion of the determinant along the columns $a_{s}, i_{k}$ of $N$,
$\det N  = \pm \det [(xI-A)$ $(\{ b_{t}, j_{k}\}, \{ a_{s}, i_{k}\})]$.
Thus, by (\ref{det N2}),
$$
\det [(xI-A)(\{ b_{t}, j_{k}\}, \{ a_{s}, i_{k}\})] = \pm c \cdot \det N[L] \cdot \det N[R] \cdot f_{k}(x)
$$
where $f_{k}(x) = \displaystyle\sum_{g=1}^{d_{k}}(-1)^{v(Q^{k}_{g})}\mbox{wt}(Q^{k}_{g}) \cdot \det N[H' \setminus Q^{k}_{g}]$.
Note that $f_{k}(x)$ is independent of the choices of $a_{s}, b_{t}$.
Since $(xI-A)(\{ b_{t}, j_{k}\}, \{ a_{s}, i_{k}\})$ is an $n-2$ by $n-2$ submatrix of $xI-A$,

$$
\Delta_{n-2}(x) | \det N[L] \cdot \det N[R] \cdot f_{k}(x).
$$

Let $\mu \in \sigma(A)$ and $m_{A}(\mu) =3$.
Then, by Theorem~\ref{factor main}, $\mu$ is a zero of
$\det N[L] \cdot \det N[R] \cdot f_{k}(x)$.
If $\mu$ is an eigenvalue of neither $A[L]$ nor $A[R]$, then
$\mu$ is a zero of $f_{k}(x)$.
Note that $L$ and $R$ are arbitrary legs of $G$ incident to different $v_{i}$'s.

If $\mu$ is not an eigenvalue of at least 4 of $A[L_{r}]$'s, then
there exist $A[L_{u}]$ and $A[L_{v}]$ such that $\lambda \not\in \sigma(A[L_{u}])$ and $\lambda \not\in \sigma(A[L_{v}])$, and
$L_{u}, L_{v}$ are connected to different $v_{i}$'s.
Thus, in this case, $\mu$ is a zero of $f_{k}(x)$ for some $k \in \{1, 2, 3\}$.
Note that $V(H') \setminus V(Q^{k}_{g}) \subseteq V(H \setminus U)$ and thereby,
the degree of $f_{k}(x)$ is $m-3$ or less for each $k =1, 2, 3$.
Therefore, the number, $n_{3}'$, of such $\mu$ satisfies
\begin{equation}
\label{n3'}
n_{3}' \leq 3(m-3).
\end{equation}

Let $n_{3}''$ be the number of eigenvalues $\mu$ of $A$ with multiplicity 3 such that $\mu$ is an eigenvalue of at least 4 of $A[L_{r}]$'s.
Then $n_{3} = n_{3}' + n_{3}''$.
Since $\mu$ is an eigenvalue of $A[L_{1}] \oplus \cdots \oplus A[L_{6}]$ with multiplicity at least 4, and the order of
$A[L_{1}] \oplus \cdots \oplus A[L_{6}]$ is $6\ell$, $4 n_{3}'' \leq 6\ell$.
Equivalently,
\begin{equation}
\label{n3''}
n_{3}'' \leq \displaystyle\frac{3\ell}{2}.
\end{equation}

Now, we compute a lower bound on $q(A)=q$.
Note that
\begin{center}
\begin{tabular}{rcl}
$n=6\ell+m = \displaystyle\sum_{i=1}^{q}{m_{A}(\lambda_{i})}$ &=& $2q -n_{1} +n_{3}' + n_{3}'' +
\displaystyle\sum_{\lambda_{j}: m_{A}(\lambda_{j}) \geq 4}{(m_{A}(\lambda_{j}) - 2)}$ \\
& $\leq$ & $2q +n_{3}' + n_{3}'' + \displaystyle\sum_{\lambda_{j}: m_{A}(\lambda_{j}) \geq 4}{(m_{A}(\lambda_{j}) - 2)}$.
\end{tabular}
\end{center}
\noindent
By (\ref{leq 2m-6}), (\ref{n3'}), and (\ref{n3''}), we have
\begin{equation}
\label{6l+m}
6\ell+m \leq 2q + 3(m -3) + \frac{3\ell}{2} + 2m -6.
\end{equation}
Hence, by solving (\ref{6l+m}) for $q=q(A)$, we have
\begin{equation}
\label{q(A)3}
q(A) \geq \frac{9\ell}{4} - 2m + \frac{15}{2}.
\end{equation}
Since $A$ is an arbitrary matrix in $S(G)$, (\ref{q(A)3}) implies

$$
q(G) \geq \frac{9\ell}{4} - 2m + \frac{15}{2}.
$$

\hfill{\rule{2mm}{2mm}}

\end{document}